\documentclass{amsart}
\usepackage{graphicx}
\usepackage{amssymb}
\usepackage{epstopdf}
\usepackage{nicefrac}
\usepackage{enumerate}
\usepackage{hyperref}
\usepackage{float}
\usepackage{color}
\usepackage{pst-all}
\usepackage{tabularx}

\newtheorem{thm}{THEOREM}[section]

\newtheorem{cor}[thm]{COROLLARY}

\newtheorem{defn}[thm]{DEFINITION}
\newtheorem{ex}[thm]{EXAMPLE}

\newtheorem{lemma}[thm]{LEMMA}
\newtheorem{prob}[thm]{PROBLEM}
\newtheorem{prop}[thm]{PROPOSITION}

\newcommand{\ds}{\displaystyle}

\newcommand{\F}{{\mathcal F}} 

 \newcommand{\bA}{{\bf A}}
 
 \newcommand{\bM}{{\bf M}}
 
\newcommand{\cG}{{\mathcal G}}
\newcommand{\cH}{{\mathcal H}}

\newcommand{\cN}{{\mathcal N}}

\newcommand{\cP}{{\mathcal P}}
\newcommand{\cQ}{{\mathcal Q}}
\newcommand{\cS}{{\mathcal S}}

\newcommand{\G}{\Gamma}
\newcommand{\mR}{{\mathbb R}}
\newcommand{\mS}{{\mathbb S}}
\newcommand{\mT}{{\mathbb T}}

\newcommand{\ove}{{\overline e}}
\newcommand{\oU}{{\overline U}}
\newcommand{\wtd}{{\widetilde d}}
\newcommand{\wtx}{{\widetilde x}}

\newcommand{\mH}{{\mathbb H}}
\newcommand{\mN}{{\mathbb N}}
\newcommand{\mZ}{{\mathbb Z}}

\newcommand{\fM}{{\mathfrak{M}}}

\newcommand{\fG}{{\mathfrak{G}}}
\newcommand{\fT}{{\mathfrak{T}}}
\newcommand{\fX}{{\mathfrak{X}}}
\newcommand{\fZ}{{\mathfrak{Z}}}
\newcommand{\wtM}{\widetilde{M}}
\newcommand{\dF}{d_{\F}} 

\input xy
\xyoption {all}

\begin{document}

\title{Growth and homogeneity of matchbox manifolds}

 \begin{abstract}   
 A matchbox manifold with one-dimensional leaves which has equicontinuous holonomy dynamics must be a homogeneous space, and so must be homeomorphic to a classical Vietoris solenoid. In this work, we consider the problem, what can be said about a matchbox manifold with equicontinuous holonomy dynamics, and all of whose leaves have at most polynomial growth type? We show that such a space must have a finite covering for which the global holonomy group of its foliation is nilpotent. As a consequence, we show that if the growth type of the leaves is polynomial of degree at most 3, then there exists a finite covering which is homogeneous. If the growth type of the leaves is polynomial of degree at least 4, then there are additional obstructions to homogeneity, which arise from the structure of nilpotent groups.

 \end{abstract}

\author{Jessica Dyer}
\author{Steven Hurder}
 \author{Olga Lukina}
 \email{jesscdyer@gmail.com, hurder@uic.edu, lukina@uic.edu}
\address{JD: Department of Mathematics, Tufts University, Bromfield-Pearson Hall, 503 Boston Avenue, Medford, MA 02155}
\address{SH \& OL: Department of Mathematics, University of Illinois at Chicago, 322 SEO (m/c 249), 851 S. Morgan Street, Chicago, IL 60607-7045}

\thanks{Version date: February 1, 2016; revised June 8, 2016}



\maketitle

 \vspace{-.1in}

 \tableofcontents

\vfill
\eject

\section{Introduction}\label{sec-intro}

 A \emph{continuum} is a compact connected metric space. A continuum $X$ is \emph{homogeneous} if for each $x, y \in X$, there exists a homeomorphism $h \colon X \to X$ such that $h(x) = y$. 
 For example, a compact connected manifold without boundary is a homogeneous continuum, and the proof of this is a standard exercise in manifold theory. On the other hand, 
 Knaster and Kuratowski \cite{KK1920} posed the problem in   1920 to classify the homogeneous continua in the plane $\mR^2$, and this problem was only recently solved by 
   Hoehn and Oversteegen   \cite{HO2014}.  Their paper also gives a   selection of references for the ``rich literature concerning homogeneous continua''. 
 For a continuum embedded in   Euclidean space $\mR^n$ with $n > 2$, the classification problem for homogeneous continua becomes intractable in this full generality, so one formulates a more restricted problem by imposing conditions on the continua.

 In this work, we are concerned with continua that are ``manifold-like''. That is,   the continuum is a disjoint union of manifolds of the same dimension, and locally has a uniform structure. This is formulated by introducing the notion of an
 \emph{$n$-dimensional foliated space} $\fM$, which  is a continuum
that has a regular  local product structure \cite{CandelConlon2000,MS2006};
that is, every point  $x \in \fM$ has an open neighborhood $x \in U \subset \fM$ homeomorphic
to an open subset of $\mR^n$ times a compact metric space $\fT_x$ where $\fT_x$ is called the local
transverse model. The homeomorphism $\varphi_x \colon \oU_x \to [-1,1]^n \times \fT_x$ is called a local foliation chart. 
A \emph{matchbox manifold} is a foliated space $\fM$ such that the
local transverse models $\fT_x$ are totally disconnected. 
The leaves of the foliation $\F$ of $\fM$ are the
maximal path connected components, and by assumption they are smooth manifolds, and can be endowed with complete Riemannian metrics.  Precise
definitions can be found in the references \cite{AM1988,CandelConlon2000,ClarkHurder2013,CHL2013c}.

 Bing conjectured  in  \cite{Bing1960} that a homogeneous continuum whose arc-composants are arcs must be a classical Van Dantzig - Vietoris solenoid \cite{vanDantzig1930,Vietoris1927}. This condition is satisfied by $1$-dimensional matchbox manifolds, so in particular Bing's   conjecture implies that such spaces   are solenoids if homogeneous.  
An affirmative answer to this  conjecture of Bing was given by Hagopian~\cite{Hagopian1977}, and
 subsequent proofs  in the framework of $1$-dimensional
matchbox manifolds were given by Mislove and Rogers \cite{MR1989} and by  Aarts,  Hagopian and Oversteegen~\cite{AHO1991}.

  Clark and Hurder   generalized  the $1$-dimensional result    to higher dimensional leaves  in the work \cite{ClarkHurder2013}, giving a positive   solution to  Conjecture~4 of      \cite{FO2002}.
 In Section~\ref{sec-solenoids}, we recall the notion of weak solenoids, and the special case of normal (or McCord) solenoids, as introduced by McCord in \cite{McCord1965}. Then the following result was proved in  \cite{ClarkHurder2013}:

\begin{thm}[Clark \& Hurder] \label{thm-homogeneous}
Let $\fM$ be a   homogeneous   matchbox manifold. Then $\fM$ is
homeomorphic to a McCord solenoid.  
\end{thm}

One step in the proof of Theorem~\ref{thm-homogeneous} is the proof of the following result, which generalizes a key step in the proofs of the $1$-dimensional case in \cite{AHO1991,MR1989}. A foliation $\F$ is said to be \emph{equicontinuous} if the transverse holonomy pseudogroup defined by the parallel transport along the leaves of $\F$ acting on a transversal space has equicontinuous dynamics. (See \cite{ClarkHurder2013} for a detailed discussion of this property.)
\begin{thm}[Clark \& Hurder \cite{ClarkHurder2013}] \label{thm-equicontinuous}
Let $\fM$ be an equicontinuous    matchbox manifold.   Then $\fM$ is
homeomorphic to a weak solenoid.  
\end{thm}
The assumption that the holonomy action is equicontinuous is essential. The Williams solenoids, as defined in \cite{Williams1974} as   inverse limits of maps between  branched $n$-manifolds,  are matchbox manifolds whose holonomy dynamics is expansive, and they are not homeomorphic to a weak solenoid.

Examples of Schori \cite{Schori1966} and Rogers and Tollefson \cite{RT1971b} show that there exists weak solenoids that are not homogeneous. The paper of Fokkink and Oversteegen \cite{FO2002} analyzed the problem of showing that a weak solenoid is in fact homogeneous, introducing a technique based on the group chains associated to a weak solenoid, to obtain a criterion for when a weak solenoid is   homogeneous.

The study of equivalence classes of group chains in the paper \cite{FO2002}  was continued in the  
   thesis \cite{Dyer2015} of the first author, and the subsequent  work \cite{DHL2016a} of the authors of this paper. As discussed in Section~\ref{sec-chains} below,     there are multiple phenomena in the algebraic structure of group chains, which may give rise to non-homogeneous weak solenoids. However, these properties of solenoids differ in their nature, and it was observed in  \cite{DHL2016a} that in some cases, the non-homogeneity property is not preserved under finite coverings, while in other cases it is. This motivated our introduction of the following notion:

    \begin{defn}\label{defn-normalcover}
 Let $\Pi_0 \colon \cS \to M_0$ be a weak solenoid with base a closed manifold $M_0$.  We say that $\cS$ is \emph{virtually homogeneous} if there  exists a finite-to-one covering map $p_0' \colon M_0' \to M_0$ such that the pullback solenoid $\Pi_0' \colon  \cS' \to M_0'$ is homogeneous.
 \end{defn}
As an example of this phenomenon, the example of Rogers and Tollefson in \cite{RT1971b} (and as interpreted by Fokkink and Oversteegen in \cite[Section~6]{FO2002}) is not homogeneous due to the    non-orientation-preserving action of an element of the fundamental group of the base manifold in their example. Passing to a finite covering removes this action, and results in a homogeneous solenoid. 

More generally, there are other properties of   weak solenoids which become evident after passing to an appropriate finite covering. In this work, we will also consider the property that a solenoid has transverse model space defined by a group chain in a nilpotent group, in which case we say that the solenoid is \emph{virtually nilpotent}.  This will be considered further in Sections~\ref{sec-holonomy} and \ref{sec-nilpotent}.

In  this work, we   consider the (virtual) homogeneity properties of   weak solenoids in terms of the geometry of their leaves.  A weak solenoid with $1$-dimensional leaves must be homogeneous, as the base manifold $M_0$ is then $1$-dimensional, hence is homeomorphic to $\mS^1$ and its fundamental group $H_0$ is free abelian.  However, already in the case of a $2$-dimensional weak solenoid, the base manifold $M_0$ is a closed $2$-manifold, possibly non-orientable. If the base is orientable, then it is either  homeomorphic to the $2$-torus $\mT^2$ or to a closed surface with genus at least $2$. Thus, if   the   leaves of the solenoid $\cS$  are simply connected, then they are either   coarse-isometric to the Euclidean plane, or to the hyperbolic plane. In the general case, a weak solenoid may have a mixture of topological types for its leaves, as  they may have differing fundamental groups for the leaves as in the examples in \cite{CFL2014, DHL2016a}, and thus cannot be homogeneous.

  The ``growth rates of   leaves'' is a standard notion in the theory of smooth foliations of manifolds \cite{Milnor1968a,Plante1975,PlanteThurston1976},  and adapts in a straightforward way to the case of matchbox manifolds. For the case of weak solenoids as discussed in section~\ref{sec-solenoids}, the growth rates of the leaves and of the fundamental group are closely related, as discussed in section~\ref{sec-growth}.   The growth rates of leaves in a weak solenoid can either be polynomial of some degree $0 \leq k < \infty$, or   can be subexponential but not polynomial, or   can be exponential.  
 Groups with polynomial growth type are special, in that their algebraic structures are understood in broad outline, as discussed in section~\ref{sec-growth}. 
This is the basis of our main results, which concern the (virtual) homogeneity properties of equicontinuous matchbox manifolds which have all leaves of polynomial growth. 
 
To state our results, we require the additional concept   of ``finite type'' for foliations. Recall that a topological space $X$ has finite type if it is homotopy equivalent to a finite CW complex. A foliated manifold has \emph{finite type} if each of its leaves is a space of finite type. This is a non-trivial assumption, as Hector constructed in  \cite{Hector1977} examples of smooth foliations of codimension-one on compact manifolds which have leaves that are not of finite type. Also,  the Schori solenoid  \cite{Schori1966} has leaves which are  surfaces of infinite genus, hence   are not of finite type.

\begin{thm}\label{thm-main0}
Let $\fM$ be an equicontinuous matchbox manifold,   and suppose that   $\F$   has finite type and all leaves of $\F$ have   polynomial growth. 
Then $\fM$  is homeomorphic to a virtually nilpotent weak solenoid. 
\end{thm}
The definitions of a virtually nilpotent and virtually abelian   solenoid are given in Definition~\ref{def-virtuallynilpotent}.

The importance of the conclusion of Theorem~\ref{thm-main0}, is that whether $\fM$ is virtually homogeneous, or not, is  then determined by algebraic invariants associated with a group chain in a nilpotent group $N_0$, and the adjoint action of a finite group $H_0$ acting on $N_0$. The  authors' works \cite{Dyer2015,DHL2016a} give a selection of examples to illustrate when such chains yield virtually homogeneous weak solenoids.

For a $1$-dimensional matchbox manifold $\fM$, all leaves have polynomial growth of degree $1$, and as noted above, equicontinuity of the holonomy of the foliation implies   that $\fM$ is homogeneous.  It is then natural to consider the applications of Theorem~\ref{thm-main0} to the virtual homogeneity   for weak solenoids,  in the cases where the leaves of $\F$ have low degree polynomial growth.  

\begin{thm}\label{thm-main1}
Let $\fM$ be an equicontinuous matchbox manifold,   and suppose that   $\F$   has finite type and all leaves of $\F$ have   polynomial growth of degree at most $3$. 
Then $\fM$  is homeomorphic to a virtually abelian weak solenoid, and thus to a virtually homogeneous weak solenoid.
\end{thm}
The proof of this result reveals that  the lack of homogeneity in the weak solenoid $\cS$  is due to the monodromy action of some finite quotient group of the fundamental group of the base of the weak solenoid, and this symmetry-breaking action can be removed by passing to a finite normal covering of the base.  
Thus, all   examples as in Theorem~\ref{thm-main1} are   analogous to the Rogers and Tollefson example.

Our third result is  an example  as constructed by Dyer in  \cite{Dyer2015}, which shows that for higher growth rates, the lack of homogeneity is due to the structural properties of the nilpotent group and the group chain in it which defines the transverse model space for the solenoid.

\begin{thm}\label{thm-main2}
There exist non-homogeneous weak solenoids with leaves of dimension $3$, where every  leaf is   covered by   the Heisenberg space $\mH$, and thus  has polynomial growth rate exactly $4$. Moreover, these examples are not virtually homogeneous.
\end{thm}
 
 These results suggest that the study of the homogeneity properties of weak solenoids with leaves of subexponential growth is itself an interesting subject for further investigation. 

 In section~\ref{sec-solenoids} we discuss the definitions and some of the properties of weak solenoids. In section~\ref{sec-growth} we discuss the relation between growth properties of leaves and the group structure for weak solenoids.  In section~\ref{sec-chains} we introduce the group chains associated to weak solenoids, and their classification as discussed in the works \cite{Dyer2015,DHL2016a,FO2002}. In  section~\ref{sec-proofs} we first prove Theorem~\ref{thm-main0}, then apply this result to obtain the proof of Theorem~\ref{thm-main1}. Finally, we give the constructions used to prove  Theorem~\ref{thm-main2}.
 The last Section~\ref{sec-future} discusses some directions for future research.
 
 The results of this paper were partially presented in the talk \cite{Hurder2015} by the second author, and also rely in part on the thesis work of the first author \cite{Dyer2015}.

\section{Weak solenoids}\label{sec-solenoids}

 In this section, we describe the constructions of  \emph{weak  solenoids}, and   recall   some of their properties.

  A  \emph{presentation}   is a collection $\cP = \{ p_{\ell+1} \colon M_{\ell+1} \to M_{\ell} \mid \ell \geq 0\}$, where each $M_{\ell}$ is a connected compact simplicial complex of dimension $n$, and each  \emph{bonding} map $p_{\ell +1}$  is a proper surjective map of   simplicial complexes with discrete fibers.
For   $\ell \geq 0$ and $x \in M_{\ell}$,  the set $\{p_{\ell +1}^{-1}(x) \} \subset M_{\ell +1}$  is compact and discrete, so   the cardinality $\# \{p_{\ell +1}^{-1}(x) \} < \infty$. For an inverse limit defined in this generality, the cardinality of the fibers  of the maps $p_{\ell +1}$  need not be constant in either  $\ell$ or $x$.  
 
    Associated to a presentation $\cP$ is an inverse limit space,  
\begin{equation}\label{eq-presentationinvlim}
\cS_{\cP} \equiv \lim_{\longleftarrow} ~ \{ p_{\ell +1} \colon M_{\ell +1} \to M_{\ell}\} ~ \subset \prod_{\ell \geq 0} ~ M_{\ell} ~ .
\end{equation}
 By definition, for a sequence $\{x_{\ell} \in M_{\ell} \mid \ell \geq 0\}$, we have 
\begin{equation}\label{eq-presentationinvlim2}
x = (x_0, x_1, \ldots ) \in \cS_{\cP}   ~ \Longleftrightarrow  ~ p_{\ell}(x_{\ell}) =  x_{\ell-1} ~ {\rm for ~ all} ~ \ell \geq 1 ~. 
\end{equation}
The set $\cS_{\cP}$ is given  the relative  topology, induced from the product topology, so that $\cS_{\cP}$ is itself compact and connected.

For example, if    $M_{\ell} = \mS^1$ for each $\ell \geq 0$, and the map $p_{\ell}$ is a proper covering map of degree $m_{\ell} > 1$ for $\ell \geq 1$, then $\cS_{\cP}$ is an example of a  {classic solenoid},  discovered independently by    van~Dantzig \cite{vanDantzig1930} and   Vietoris   \cite{Vietoris1927}.
If     $M_{\ell}$ is a compact manifold without boundary for each $\ell \geq 0$,  and  the map $p_{\ell}$ is a proper covering map of degree $m_{\ell} > 1$  for   $\ell \geq 1$,   then $\cS_{\cP}$  is said to be a \emph{weak solenoid}. This generalization of $1$-dimensional solenoids was originally    considered  in the papers by McCord \cite{McCord1965} and Schori  \cite{Schori1966}. In particular, McCord showed in \cite{McCord1965} that   $\cS_{\cP}$ has a local product structure.

\begin{prop}\label{prop-solenoidsMM}
Let   $\cS_{\cP}$ be a weak solenoid, whose    base space $M_0$   is a compact manifold of dimension $n \geq 1$. Then   $\cS_{\cP}$ is  a minimal matchbox manifold of dimension $n$ with foliation   $\F_{\cP}$. 
\end{prop}

Associated to a presentation $\cP$ of compact manifolds is a sequence of proper surjective maps 
\begin{equation}\label{eq-coverings}
q_{\ell} = p_{1} \circ \cdots \circ p_{\ell -1} \circ p_{\ell} \colon M_{\ell} \to M_0 ~ .
\end{equation}
For each $\ell > 1$, projection onto the $\ell$-th factor in the product $\ds \prod_{\ell \geq 0} ~ M_{\ell}$ in \eqref{eq-presentationinvlim} yields a 
  fibration map denoted by $\Pi_{\ell} \colon \cS_{\cP}  \to M_{\ell}$, for which 
 $\Pi_0 = \Pi_{\ell} \circ q_{\ell} \colon \cS_{\cP} \to M_0$. 
A choice of a basepoint $x_0 \in M_0$ fixes a fiber  $\fX_0 = \Pi_0^{-1}(x_0)$, which is a Cantor set by the assumption  that the fibers of each map $p_{\ell}$ have cardinality at least $2$.  The choice of $x_0$ will remain fixed throughout the following. We also then have a fixed group   $H_0 =  \pi_1(M_{0}, x_{0})$.

A choice  $x \in \fX_0$   defines basepoints  $x_{\ell} = \Pi_{\ell}(x) \in M_{\ell}$ for $\ell \geq 1$.
Define $\cH^x_{\ell} = \pi_1(M_{\ell}, x_{\ell})$, and let
\begin{equation}\label{eq-imahes}
H^x_{\ell} = {\rm image}\left\{  (q_{\ell} )_{\#} \colon \cH^x_{\ell} \longrightarrow H_{0}\right\}  
\end{equation}
  denote  the image of the induced map $(q_{\ell} )_{\#} $ on fundamental groups. Thus, associated to the presentation $\cP$ and basepoint $x \in \fX_0$ we obtain a descending chain of subgroups of finite index
  \begin{equation}\label{eq-descendingchain}
H_{0} \supset H^x_{1} \supset H^x_{2} \supset \cdots \supset H^x_{\ell} \supset \cdots  \ .
\end{equation}
Each quotient  $X_{\ell}^x = H_{0}/H_{\ell}^x$ is   a   finite set equipped with a left $H_0$-action, and there are surjections $X_{\ell +1}^x \to X_{\ell}^x$ which commute with the action of $H_0$.  The inverse limit 
\begin{equation}\label{eq-Galoisfiber}
X_{\infty}^x = \lim_{\longleftarrow} ~ \{ p_{\ell +1} \colon X_{\ell +1}^x \to X_{\ell}^x\} ~ \subset \prod_{\ell \geq 0} ~ X_{\ell}^x  
\end{equation}
is then a totally disconnected perfect set, so is   a Cantor set. The   fundamental group $H_0$ acts on the left on  $X_{\infty}^x$ via     the coordinate-wise multiplication on the product in \eqref{eq-Galoisfiber}.    We denote this   Cantor action by $(X_{\infty}^x , H_0 , \Phi_x)$.

 \begin{lemma}\label{lem-denseaction}
The left action   $\Phi_x \colon H_0 \times X_{\infty}^x \to X_{\infty}^x$ is minimal. 
\end{lemma}  
 \proof
 The left action of $H_0$ on each quotient space $X_{\ell}^x$ is transitive, so the orbits are dense in the product topology on $X_{\infty}^x$.
  \endproof

The choice of the basepoint $x \in \cS_{\cP}$ defines   basepoints $x_{\ell} \in M_{\ell}$ for all $\ell \geq 1$, which gives   an identification of $X_{\ell}^x$ with the fiber of the covering map $M_{\ell} \to M_0$. In the inverse limit, we thus obtain a map $\tau_x \colon X_{\infty}^x \to  \fX_0 = \Pi_0^{-1}(x_0)$ which is a homeomorphism. The left action of $H_0$ on $X_{\infty}^x$ is conjugated to an action of $H_0$ on $\fX_0$  called the \emph{monodromy action} at $x_0$ for the fibration $\Pi_0 \colon \cS_{\cP} \to M_0$.   The monodromy action can also be  defined by the holonomy transport along the leaves of the foliation $\F_{\cP}$ on $\cS_{\cP}$. It was shown in \cite[Theorem~4.8]{CHL2013c} that the monodromy actions for two homeomorphic solenoids are return equivalent, where return equivalence  is a notion of Morita equivalence for matchbox manifolds.

Note that while the group chain   in \eqref{eq-descendingchain}   depends on the choice of basepoint $x$,  for $x \ne y \in \fX_0$  the composition $\tau_y^{-1} \circ \tau_x \colon X_{\infty}^x \to X_{\infty}^y$ gives a topological conjugacy between the minimal Cantor actions  $(X_{\infty}^x , H_0 , \Phi_x)$ and $(X_{\infty}^y , H_0 , \Phi_y)$.
That is, the map $\tau_x \colon X_{\infty}^x \to  \fX_0$ can be viewed as ``coordinates'' on the inverse limit space $\fX_0$, and so properties of the minimal Cantor action 
$(X_{\infty}^x , H_0 , \Phi_x)$ which are independent of the choice of coordinates are properties of the topological type of $\cS_{\cP}$. In Section~\ref{sec-chains}, we discuss   properties  of $(X_{\infty}^x , H_0 , \Phi_x)$ derived from the group chain \eqref{eq-descendingchain} and their dependence on the  choices made which are used to define it.

Let $\wtM_0$ denote the universal covering of the compact manifold $M_0$ and let $(X_{\infty}^x , H_0 , \Phi_x)$ be the minimal Cantor action associated to the presentation $\cP$ and the choice of a basepoint $x \in \fX_0$. Associated to the left action $\Phi_x$ of $H_0$ on $X_{\infty}^x$ is a suspension space 
\begin{equation}\label{eq-suspensionfols}
\fM = \wtM_0 \times X_{\infty}^x / (z \cdot g^{-1}, x) \sim (z , \Phi_x(g)(x)) \quad {\rm for }~ z \in \wtM_0 , ~ g \in H_0 ~,
\end{equation}
which is a minimal matchbox manifold. 
Given  coverings $\pi' \colon M' \to M_0$ and $\pi'' \colon M'' \to M_0$, and choices of basepoints $x' \in M'$ and $x'' \in M''$ with $\pi'(x') = \pi''(x'') = x_0$, such that    the subgroups 
$$\ds \pi_{\#}'(\pi_1(M',x')) = \pi_{\#}''(\pi_1(M'', x'')) \subset \pi_1(M_0,x_0) ,$$ 
then there is a natural homeomorphism of coverings $M' \cong M''$ which is defined using the path lifting property. From this, it   follows (see   \cite{ClarkHurder2013}) that:
\begin{thm}\label{thm-weaksuspensions}
Let   $\cS_{\cP}$ be a weak solenoid with base space $M_0$ where $M_0$ is a compact manifold of dimension $n \geq 1$. Then there is a foliated homeomorphism $\cS_{\cP} \cong \fM$.
\end{thm}
 
 \begin{cor}\label{cor-weaksuspensions}
The homeomorphism type of a weak solenoid    $\cS_{\cP}$ is completely determined by the base manifold $M_0$ and the associated Cantor action by $(X_{\infty}^x , H_0 , \Phi_x)$.
\end{cor}

Note that all leaves in a McCord solenoid $\cS_{\cP}$ are homeomorphic, as it is homogeneous.  
In the case of a weak solenoid $\cS_{\cP}$, the leaves of $\F_{\cP}$ need not be homeomorphic. For example, the work \cite{CFL2014} gives an example of a weak solenoid  for which   the leaves of $\F_{\cP}$   have differing numbers of ends.

\section{Growth properties}\label{sec-growth}

In this section, we discuss the growth properties of leaves of foliated spaces and of finitely generated groups, and the relations between these two concepts in the case of weak solenoids. 

Let $\fM$ be a matchbox manifold. 
A map $f \colon \fM \to \mR$ is said to be \emph{smooth} if for each flow box
$\varphi_x \colon \oU_x \to [-1,1]^n \times \fT_x$ for $\fM$ and $w \in \fT_x$ the composition
$y \mapsto f \circ \varphi_x^{-1}(y, w)$ is a smooth function of $y \in (-1,1)^n$, and depends continuously on $w$ in the $C^{\infty}$-topology on maps of the plaque coordinates $y$. As noted in \cite{MS2006}, and also \cite[Chapter 11]{CandelConlon2000}, this allows one to define smooth partitions of unity, vector bundles, and tensors for smooth foliated spaces. In particular, one can define leafwise Riemannian metrics, which are defined for vectors in the tangent bundle $T\F$ to $\F$, and so define a norm on the tangent space $T_x\F$ for $x \in \fM$. This norm is smooth as $x$ varies in a leaf, and is continuous for $x$ in $\fM$.
We then recall a standard result, whose basic idea dates back to the work of Plante \cite{Plante1975} if not before. The proof for foliated spaces can be found in \cite[Theorem~11.4.3]{CandelConlon2000}.
\begin{thm}\label{thm-riemannian}
Let $\fM$ be a smooth foliated space. Then there exists a leafwise Riemannian metric for $\F$, such that for each $x \in \fM$, the leaf $L_x$ through $x$ is a complete Riemannian manifold with bounded geometry, and the Riemannian geometry depends continuously on $x$ . 
\end{thm}

The survey \cite{Hurder1994} discusses  results and  problems concerned with the ``coarse geometry'' of smooth foliations, and many of these results apply as well to the study of foliated spaces. We recall some basic concepts of coarse geometry. 

Let $(X, d_X)$ and $(Y, d_Y)$ be metric spaces. A homeomorphism $f \colon X \to Y$ is said to be  \emph{bi-Lipschitz} if there exists a constant $C \geq 1$ such that for all $x,x' \in X$ we have
\begin{equation}\label{eq-Lipschitz}
C^{-1} \cdot d_X(x, x') \leq d_Y(f(x),f(x')) \leq C \cdot d_X(x, x') \  .
\end{equation}
The map $f$ is said to be an \emph{isometry} if the condition \eqref{eq-Lipschitz} is true for $C=1$. 

A \emph{set map} $f \colon X \to Y$ is said to be a \emph{quasi-isometry}, or more precisely a $(C,D)$-quasi-isometry,  if there exists  constants $C \geq 1$ and $D \geq 0$ such that for all $x,x' \in X$ we have
$$ C^{-1} \cdot d_X(x, x') - D \leq d_Y(f(x),(x')) \leq C \cdot d_X(x, x') + D$$
and the image of $f$ is $D$-dense in $Y$. Recall that a subset $Z \subset Y$ is said to be $D$-dense if for all $y \in Y$, there exists $z \in Z$ such that $d_Y(y, z) \leq D$.
In the case of the image of a   map $f$, this means that  for all $y \in Y$, there exists $x \in X$ such that $d_Y(y, f(x)) \leq D$.

The coarse geometry of  metric spaces is the study of the properties of the metric geometry which are preserved by a quasi-isometry. One of the basic ideas is the notion of a \emph{net}, or \emph{Delone set} as it is called in the tilings literature. Let $(X, d_X)$ be a metric space. Then a subset $Z \subset X$ is a net if there exists constants $A,B > 0$ such that: 
\begin{itemize}
\item For all $x \in X$ there exists $z \in Z$ such that $d_X(x,z) \leq A$.  
\item For all $z \ne z' \in Z$ then $d_X(z,z') \geq B$.
\end{itemize}
It follows that for the   metric $d_Z$ on $Z$ obtained  from the restriction of $d_X$, the space $(Z, d_Z)$ is discrete, and the inclusion $Z \subset X$ is a $(C,D)$-coarse isometry where $C = 1$ and $D = A$. 

The simplest example of an $A$-dense set is the inclusion $\mZ \subset \mR$ for the standard metric on $\mR$, where $A = 1/2$. Moreover, the subset $\mZ$ is a net in $\mR$.
This is a special case of a more general class of examples. Let $M$ be a closed Riemannian manifold, with the path metric $d_M$. Let $\pi \colon \wtM \to M$ denote its universal covering space. Then the Riemannian metric on $M$ lifts to a Riemannian metric on $\wtM$, and   the induced path length metric $\wtd$ on $\wtM$ is complete. Let $x \in M$ and $Z_x = \pi^{-1}(x) \subset \wtM$ be the fiber over $x$. Then $Z_x$ is a net for $(\wtM, \wtd)$ for $A = {\rm diam}(M)/2$ and $B = {\rm inj}(M)$, where ${\rm diam}(M)$ is the diameter of $M$ for the metric $d_M$ and ${\rm inj}(M)$ in the injectivity radius for $M$, which is positive as $M$ is compact without boundary.   This was a key observation in the seminal work by Milnor \cite{Milnor1968a}. For the covering $\pi \colon \mR \to \mS^1$ with the usual metric on $\mS^1$, this construction yields the example above of the inclusion $\mZ \subset \mR$.

Another class of examples of coarse geometry are provided by the leaves of a foliation of a manifold, or of a foliated space $\fM$. Assume that a leafwise Riemannian metric has been chosen for $\fM$. For each leaf $L \subset \fM$, let $d_L$ denote the path length metric on $L$ defined by the restriction of the Riemannian metric to $L$. Let $\dF$ denote the metric on $\fM$ defined by the collection of leafwise metrics for the leaves. That is, if $x, y \in \fM$ and there is some leaf $L$ with $x,y \in L$ then $\dF(x,y) = d_L(x,y)$. Otherwise, we set $\dF(x,y) = \infty$. We call $\dF$ a \emph{leafwise metric} for $\fM$.   

  Plante observed in \cite{Plante1975} that for smooth foliations of compact manifolds,    for any two choices of leafwise Riemannian metrics on $\F$, the resulting leafwise metrics $\dF$ and $\dF'$ on $M$ are bi-Lipschitz equivalent, where the Lipschitz constant $C \geq 1$ depends on the choice of the metrics.  This remains true for leafwise metrics on foliated spaces and matchbox manifolds.

If $\fZ \subset \fM$ is a complete transversal to $\F$, then for each leaf $L \subset \fM$ the intersection $Z = \fZ \cap L$ is a net for $(L, d_L)$, so the inclusion map $\fZ \subset \fM$ induces a coarse isometry when restricted to leaves.
From this it follows that if $\fM$ is a matchbox manifold   with leafwise metric $\dF$, and $\fM'$  is a matchbox manifold with leafwise metric $\dF'$, then a homeomorphism $h \colon \fM \to \fM'$ induces a coarse isometry between each leaf $L \subset \fM$ and the leaf $L' \subset \fM'$ for which $h(L) \subset L'$. We conclude:
\begin{thm}\label{thm-coarseinvariance}
Let $h \colon \fM \to \fM'$ be a homeomorphism between matchbox manifolds. Let $L \subset \fM$ be a leaf, and let $L' = h(L) \subset \fM'$ be the image leaf. Then for any choice of leafwise Riemannian metrics on $\fM$ and $\fM'$, the map $h \colon L \to L'$ induces a coarse isometry between the leaves.
\end{thm}

Theorem \ref{thm-coarseinvariance} states that the coarse geometry of leaves in a matchbox manifold are preserved by homeomorphism of the ambient foliated spaces. It thus makes sense to consider invariants of the coarse geometry of the leaves of a matchbox manifold $\fM$, and the relations between  these invariants and the   group ${Homeo}(\fM)$.

In the case of a weak solenoid, there is a standard construction of a leafwise metric $\dF$ on the matchbox manifold $\fM =  \cS_{\cP}$. Let $\cS_{\cP}$ be   defined by a presentation $\cP = \{ p_{\ell+1} \colon M_{\ell+1} \to M_{\ell} \mid \ell \geq 0\}$, with 
 $\Pi_{0} \colon \cS_{\cP}  \to M_{0}$ the projection map onto the base manifold $M_0$. Recall that $\F_{\cP}$ denotes the resulting foliation on $\cS_{\cP}$. 
Choose a Riemannian metric on $TM_0$. Then for each leaf $L \subset \cS_\cP$ of $\F_{\cP}$, the restriction $\Pi_0 \colon L \to M_0$ is a covering map, so defines a smooth Riemannian metric on the tangent bundle $TL$. Let $d_L$ denote the path metric induced on $L$, which is complete as $M_0$ has no boundary. 
The collection of metrics so defined yields a continuous Riemannian metric on the tangential distribution $T\F_{\cP}$ to   $\F_{\cP}$, and  $\dF$ is the resulting leafwise metric.

We next introduce the \emph{growth type}  of a complete Riemannian manifold $L$, with $d_L$ the associated path length metric.   Given $x \in L$ and $r > 0$ let 
$$B_L(x,r) = \{ y \in L \mid d_L(x,y) \leq r\}$$
denote the ball in $L$ of radius $r$ centered at $x$. The Riemannian metric on $L$ defines a volume form, and let $Vol(B(x,r))$ denote the volume of the ball for this volume form. 
 The growth function of $L$ at $x$ is defined to be the function ${\rm Gr}(L,d_L, x,r) = Vol(B(x,r))$. This function depends on   the choices made, so we introduce an equivalence relation on     functions, which is used to define the growth type of $L$ which is a quasi-isometry invariant.
 
     Given     functions $f_1, f_2 \colon [0,\infty) \to [0, \infty)$ say that $f_1 \lesssim f_2$ if there exist constants $A, B, C > 0$ such that for all $r \geq 0$, we have that $  f_2(r) ~ \leq ~ A \cdot f_1(B \cdot r) + C$.
 Say that    $f_1 \sim f_2$ if both $f_1 \lesssim f_2$ and $f_2 \lesssim f_1$ hold.    This defines  an equivalence relation on functions, which defines their   \emph{growth class}. 
 The same definitions can also be applied to functions defined on any subset $\Lambda \subset \mR$.   In the following we will consider functions defined on the natural numbers   $\mN \subset \mR$, for example, and speak of their growth class.   
   
  One can consider a variety of special classes of growth types. For example, note that if $f_1$ is the constant function and $f_2 \sim f_1$ then $f_2$ is constant also. 
  
We say that $f \colon [0,\infty) \to [0, \infty)$ has \emph{exponential growth type} if $f(r) \sim \exp(r)$. Note that $\exp(\lambda \cdot r) \sim \exp(r)$ for any $\lambda > 0$, so there is only one growth class of ``exponential type''.
  
  A function $f$ has \emph{nonexponential growth type} if $f(r)  \lesssim \exp(r)$,  but $\exp(r) \not\lesssim f(r)$. 
  
  A function $f$ has   \emph{subexponential growth type},  if for any $\lambda > 0$ there exists $A, C > 0$ so that $f(r) \leq A \cdot \exp(\lambda \cdot r) + C$.

  Finally, $f$ has  \emph{polynomial growth type} if there exists $d \geq 0$ such that $f(r)  \lesssim r^d$. The growth type is exactly polynomial of degree $d$ if $f(r) \sim r^d$. 
  
From the remarks above in this section, we then conclude the following standard result, which is proved for foliated manifolds by Plante \cite{Plante1975}, and see also \cite{Hurder1994}. The extension to leaves of foliated manifolds and matchbox manifolds is an exercise.
\begin{prop}\label{prop-growthtype}
 Let $(L, d_L)$ be a complete Riemannian manifold. Then the growth class of the  function ${\rm Gr}(L,d_L, x,r)$ is well-defined up to quasi-isometry. In particular, if $L \subset \fM$ is a leaf of a matchbox manifold, then there is a well-defined  growth class for $L$ associated to a leafwise metric $\dF$ on $\fM$, which is independent of the choice of Riemannian metric on leaves, and basepoint $x \in L$.
\end{prop}
  
  The question we address in this work, is how are the properties of the homeomorphism group ${Homeo}(\fM)$ related to the growth classes of the leaves of   a weak solenoid? 
  To address this question, we recall the notions of growth class for finitely generated groups and nets, and their relation to the growth class of leaves. Chapters VI and VII   of the book \cite{delaHarpe2000} 
contain a concise overview of the geometric theory of finitely generated groups.

  Let $\G$ be a finitely generated group, and let $\G^{(1)} = \{ \gamma_0 = 1, \gamma_1, \ldots , \gamma_k\}$ be a set of generators. Then $\gamma \in \G$ has word norm  $\|\gamma \| \leq \ell$  (with respect to the set $\G^{(1)}$) 
  if we can express $\gamma$ as a product of at most $\ell$ generators, 
 $\gamma = \gamma_{i_1}^{\pm} \cdots \gamma_{i_{\ell}}^{\pm}$. Define the ball of radius $\ell$ about the identity of $\G$ by 
   $$\G^{(\ell)}  \equiv \{\gamma \in \G \mid \|\gamma\| \leq \ell\}.$$
   Note that $\G$ finitely generated implies that the cardinality of the set $\G^{(\ell)}$ is finite for all $\ell \geq 1$, so we can define the 
  growth function  $Gr(\G, \G^{(1)}, \ell) = \# \G^{(\ell)}$. This function depends upon the choice of generating set for $\G$, but its growth class does not, due to the following elementary result.

\begin{lemma}[\cite{Milnor1968a}]
Let $\G^{(1)} = \{\gamma_1, \ldots , \gamma_k\}$ and $\G^{(1)'} = \{\gamma_1', \ldots , \gamma_{k'}'\}$ be two sets of generators for $\G$. Let $B\geq 1$ be an integer such that each generator $\gamma_j'$ can be expressed by a word in at most $B$ elements of $\{\gamma_1^\pm, \ldots , \gamma_k^\pm\}$. Then 
$\ds Gr(\G', \G^{(1)'}, \ell) \leq Gr(\G, \G^{(1)}, B \cdot \ell)$.
\end{lemma}
For a finitely generated group $\G$ it thus makes sense to speak of its growth class, denoted by $Gr(\G)$. There is also the following standard result, whose proof is again elementary:
\begin{lemma}\label{lem-growtheq}
Let $H \subset \G$ be a subgroup of finite index. Then   $Gr(H) = Gr(\G)$.
\end{lemma}

Next, let    $H \subset \G$ be an arbitrary  subgroup, and let $\Lambda = \G/H$ be the     set of cosets of $H$, with the left $\G$ action. Let $\ove \in \Lambda$ denote the coset of the identity element.  The choice $\G^{(1)}$ of a generating set of $\G$ defines a word metric on $\G$ as above, which can then be used to define a quotient metric on $\Lambda$, as follows: given $x,y \in \Lambda$ set 
$$d_{\Lambda}(x,y) = \inf \, \left\{\| \gamma\| \mid \gamma \cdot x = y  \right\} .$$
If $H$ is a normal subgroup of $\G$, then $\Lambda$ is a group, and the set of generators for $\G$ projects to a set of generators for $\Lambda$. Then $d_{\Lambda}$ is simply the word metric on $\Lambda$ for this projected metric. For the case when $H$ is not a normal subgroup, then $d_{\Lambda}$ is the \emph{Cayley metric} for the discrete space $\Lambda$. We  define the growth function of the set $(\Lambda, d_{\Lambda})$ by counting the number of points in a ball of radius $\ell$, 
$$Gr(\Lambda, \G^{(1)}, \ell) = \#  \, \left\{ x \in \Lambda \mid  d_{\Lambda}(\ove, x) \leq \ell \,  \right\}$$
Let $Gr(\Lambda)$ denote the growth class of the function $Gr(\Lambda, \G^{(1)}, \ell)$, which as before, is independent of the choice of generating set $\G^{(1)}$ for $\G$. The following is then immediate:
\begin{lemma}\label{lem-growthineq}
Let $H \subset \G$ be a subgroup, and $\Lambda = \G/H$. Then   $Gr(\Lambda) \leq  Gr(\G)$.
\end{lemma}

 One of the main problems in geometric group theory is to determine the algebraic properties of a finitely generated group $\G$ which depend only on the quasi-isometry class of  $\G$ for the word metric as defined above. 
The following is a  celebrated  theorem  of Gromov.
\begin{thm} [\cite{Gromov1981}] \label{thm-Gromov}
Suppose $\G$ has polynomial growth type for some generating set. Then there exists  a subgroup of finite index $\G_0 \subset \G$ such that $\G_0$ is a nilpotent group.
  \end{thm}

 Given a nilpotent group $\G_0$, the lower central series for $\G_0$ is defined as follows. 
Let $\G_1 = [\G_0 , \G_0] \subset \G_0$ be the first commutator subgroup, which is characterized as the minimal normal subgroup $H \subset \G_0$ such that $\G_0/H$ is abelian. 
Let $r_0 \geq 0$ denote the rank of the free abelian summand of $\G_0/\G_1$. 

Then for $\ell \geq 1$,  recursively define  $\G_{\ell +1} \subset \G_{\ell}$ as the subgroup generated by all commutators:
$$[\G_0 , \G_{\ell}] = \{g^{-1}h^{-1}gh \mid  g \in \G_0 , h \in \G_{\ell} \} \ .  $$

Then $\G_{\ell +1}$ is again a normal subgroup of $\G_0$ and each quotient group $\G_{\ell}/\G_{\ell+1}$ is abelian. Let $r_{\ell} \geq 0$ denote the rank   of the free abelian summand of $\G_{\ell}/\G_{\ell + 1}$.  
 The assumption that $\G_0$ is nilpotent implies that there is a finite least index $k \geq 1$ such that $\G_{k-1}$ is non-trivial, while $\G_{k}$ is the trivial group. This is called the   \emph{length} $k(\G_0)$ of $\G_0$. If the length $k(\G_0)=1$ then $\G_0$ is abelian.  
 
 Note that if  $r_{\ell}=0$ then $\G_{\ell}/\G_{\ell+1}$ is a finitely generated torsion abelian group, but need not be the trivial group.  We recall two elementary observations.
 
   \begin{lemma}\label{lem-torsion}
 Let $\G_0$ be a finitely-generated nilpotent group,  and suppose that $r_{\ell}=0$ for $\ell \geq \ell_0$. Then $\G_{\ell_0} \subset \G_0$ is a finite nilpotent normal subgroup.
 \end{lemma}
\proof Suppose that $A \to B \to C$ is an exact sequence of groups,   with $A$ and $C$ finite,   then $B$ is also a finite group. Let $k = k(\G_0)$ then 
$\G_k = \{0\}$ and so $\G_{k-1} = \G_{k-1}/\G_{k}$ is a finite group if $k-1\geq \ell_0$. Next suppose that $k-2 \geq \ell_0$, then 
  we have the exact sequence $\G_{k-1}   \to \G_{k-2} \to \G_{k-2}/\G_{k-1}$ so that   $\G_{k-2}$ is a finite group. We continue recursively in this manner with the exact sequences $\G_{\ell}   \to \G_{\ell-1} \to \G_{\ell-1}/\G_{\ell}$ and the claim follows.
\endproof

 We also need a simple observation about the   ranks $r_0$  and   $r_1$.
  \begin{lemma}\label{lem-rankest}
 If $r_{1} \geq 1$   then $r_{0} \geq 2$.
 \end{lemma}
 \proof
 Note the commutator identity $[x,yz] = [x,y][y,[x,z]][x,z]$   for elements $x,y,z \in H$, for any group $H$. 
 In the case $H = \G_0/\G_2$ all double commutators are the identity, so this by induction yields $[x,y^{p}] = [x,y]^{p}$ for all $p > 1$. Thus, if $y^{p} = id$ for some $p> 1$, then $[x,y]^{p} = e$ also, where $e \in \G_0/\G_2$ denotes the identity element.

 Suppose $r_0 = 0$. Then for every $y \in \G_0$ there is $p>0$ such that $y^p$ is the product of commutators. Then, using the identity above, $[x,y^p]$ is the product of double commutators, and so equals $e$ in $\G_0/\G_2$. Since $[x,y^p] = [x,y]^p$ in $\G_0/\G_2$, then every commutator in $\G_0/\G_2$ is torsion. Thus,   $\G_1/\G_2$ is generated by torsion elements, and as $\G_1/\G_2$ is abelian, it is a torsion group. This contradicts the assumption that $r_1 = 1$.

Suppose $r_0 = 1$, and let $y \in \G_0$ be a generator of the free abelian factor in $\G_0/\G_1$. Then for every element $x \in \G_0$ there is a power $s >0$ such that $x^s$ is the product of $y^p$ and a finite number of commutators, where $p \geq 0$. Then we calculate in $\G_0/\G_2$ that 
$$[y^k,x^s] = [y^k,y^p][y^k,[.,.]] .. [y^k,[..,..]] = e \ , $$
and so $[y^k,x]^s = e$Ê in $\G_1/\G_2$. If $z \notin y^k\G_1$ for some $k \ne 0$, then there is $t >0$ such that $z^t$ is the product of commutators, and so $[x^s,z^t] =[x^s,z]^t= e$ in $\G_0/\G_2$. It follows that every commutator is a torsion element in $\G_0/\G_2$, and so $\G_1/\G_2$ is torsion, which contradicts the assumption that $r_1 = 1$.
\endproof

Finally, we recall a basic result of  Guivarc'h \cite{Guivarch1970,Guivarch1971} and  Bass  \cite{Bass1972} about   growth types. 

\begin{thm}[\cite{Bass1972,Guivarch1971}]\label{thm-Bass}
Let $\G_0$ be a finitely generated nilpotent group. Then the growth type of $\G_0$ is polynomial with degree 
\begin{equation}\label{eq-growthformula}
d(\G_0) = \sum_{\ell = 0}^{k(\G_0)} ~ (\ell +1)  \cdot r_{\ell} \ .
\end{equation}
\end{thm}

We apply   formula \eqref{eq-growthformula} for the cases where $d(\G_0)$ is small, as   it then strongly proscribes the ranks $r_{\ell}$.
\begin{prop}\label{prop-lowrank}
Let $\G_0$ be a finitely generated nilpotent group. If $1 \leq d(\G_0) \leq 3$, then there exists a finite-index, free abelian subgroup $\G_0' \subset \G_0$ with rank $d(\G_0)$.  
\end{prop}
\proof 
For $d(\G_0) \leq 3$,   Lemma~\ref{lem-rankest} and  formula \eqref{eq-growthformula} imply that $r_{\ell} =0$ for $\ell \geq 1$, and thus  $1 \leq r_0 = d(\G_0) \leq 3$. Moreover,    for each $\ell \geq 1$,  the quotient group $\G_{\ell}/\G_{\ell +1}$ is a finite torsion abelian group, and thus  $\G_1$ is a finite torsion group by Lemma~\ref{lem-torsion}.  

On the other hand, $r_0 \geq 1$ implies that $\G_{0}/\G_{1}$ contains a free abelian subgroup of rank $r_0$. We consider the case $r_0 = 3$, and the other cases are similar. Choose $\{\gamma_1, \gamma_2, \gamma_3\} \subset \G_0$ which map onto the generators of this free abelian subgroup of $\G_{0}/\G_{1}$. Then each commutator $[\gamma_i, \gamma_j] \in \G_1$ for $1 \leq i < j \leq 3$ is a torsion element.
Therefore, there exists an integer $m > 0$ such that each commutator $[\gamma_i^m, \gamma_j^m]$ is the trivial element. Let $\G_0'$ be the subgroup of $\G_0$ generated by $\{\gamma_1^m , \gamma_2^m, \gamma_3^m\}$. Then $\G_0'$ is free abelian with rank $3$, and has finite index in $\G_{0}/\G_{1}$ hence has finite index in $\G_0$. 
\endproof

Finally, we discuss the properties of the growth type of a covering of a compact manifold $M_0$. Let $\G_0 = \pi_1(M_0, x_0)$ for a basepoint $x_0 \in M_0$. Given a subgroup $\G_0' \subset \G_0$ let $\pi \colon \wtM_0' \to M_0$ be the associated covering space. Then for $\wtx \in \wtM_0'$ with $\pi(\wtx) = x_0$   we have 
$$\G_0 = {\rm image}\left\{ \pi_{\#} \colon \pi_1(\wtM_0' , \wtx) \to \pi_1(M_0 , x_0) \right\} .$$
The basic observation is that $M_0$ compact implies that the subset $\Lambda = \pi^{-1}(x_0) \subset \wtM_0'$ is a net,  so that with the induced metric on $\Lambda$, the growth type of $\wtM_0'$ as a Riemannian manifold  and the growth type of $\Lambda$ as a discrete subspace of $\wtM_0'$   agree. 
Moreover, covering space theory identifies $\Lambda = \G_0/\G_0'$. 
We then have the fundamental result  of Milnor and   \v{S}varc:
\begin{thm}[\cite{Milnor1968a,Svarc1955}]\label{thm-Milnor}
Let $M_0$ be a compact Riemannian manifold with basepoint $x_0 \in M_0$ and fundamental group $\G_0 = \pi_1(M_0, x_0)$. Let $\G_0' \subset \G_0$ be a   subgroup and $\pi \colon \wtM_0' \to M_0$ be the associated covering space with pathlength metric $d_{\wtM}$ on $\wtM_0'$.  Then for the quotient space  $\Lambda= \G_0/\G_0'$,  the growth type of $\wtM_0'$  equals the growth type of $\Lambda$ as a coset space of $\G_0$. 
\end{thm}

By combining  Theorem~\ref{thm-Gromov} with Theorem~\ref{thm-Milnor} for the case $\G_0'$ the trivial subgroup, we obtain:

\begin{cor} Let $\fM = \cS_{\cP}$ be a weak solenoid with base manifold $M_0$ which has a     simply connected leaf with polynomial growth for some leafwise metric. Then the fundamental group $\G_0 = \pi_1(M_0, x_0)$ has a nilpotent subgroup $\G_0' \subset \G_0$ of finite index.
\end{cor}

\section{Group chains}\label{sec-chains}

 In this section, we consider  group chains and some of their properties, and discuss the applications to the study of the automorphism groups of equicontinuous Cantor minimal systems.

\begin{defn}\label{gr-chain}
Let $G$ be a finitely generated group. A \emph{group chain} $\cG = \{G_i\}_{i \geq 0}$, with $G_0=G$, is a properly descending chain of subgroups of $G$, such that $|G:G_i| < \infty$ for every $i \geq 0$.
Let $\fG$  denote the collection of all possible nested group chains for $G$.
\end{defn}
Associated to a group chain $\{G_i\}_{i \geq 0}$, there exists a Cantor space 
  \begin{align*}   G_\infty = \lim_{\longleftarrow}\{G/G_i \to G/G_{i-1}\} = \{(G_0, g_1G_1, g_2G_2,\ldots) \mid g_jG_i = g_iG_i \textrm{ for all }j \geq i\}. 
  \end{align*}  
There is  a  natural basepoint $(eG_i) \in G_\infty$ corresponding to the  identity element $e \in G$, and a
natural left $G$-action $\Phi$ on $G_\infty$, given by 
   \begin{align}\label{inv-limaction}   \Phi(\gamma)(g_iG_i)  = \gamma \cdot (g_iG_i) =   (\gamma g_i G_i), \textrm{ for all } \gamma \in G .
   \end{align}
The $G$-action is transitive on each finite factor $G/G_i$, which implies that the $G$-action is minimal on $G_\infty$.  Moreover, the action is equicontinuous. Thus, given a group chain   $\{G_i\}_{i \geq 0}$ we obtain an equicontinuous Cantor minimal system, denoted by $(G_\infty, G, \Phi)$.

Given the group $G$, recall that $\fG$ denotes the collection  of all possible group chains in $G$.
Then there are two equivalence relations   defined on $\fG$. The first   was defined by Rogers and Tollefson in \cite{RT1971b}, and used by Fokkink and Oversteegen in \cite{FO2002} in a slightly different form.
\begin{defn} \cite{RT1971b}\label{defn-greq}
In a finitely generated group $G$, two group chains $\{G_i\}_{i \geq 0}$ and $\{H_i\}_{i \geq 0}$ with $G_0=H_0=G$ are \emph{equivalent}, if and only if, there is a group chain $\{K_i\}_{i \geq 0}$ and infinite subsequences $\{G_{i_k}\}_{k \geq 0}$ and $\{H_{j_k}\}_{k \geq 0}$ such that $K_{2k} = G_{i_k}$ and $K_{2k+1} = H_{j_k}$ for $k \geq 0$.
\end{defn}

The next definition was given by Fokkink and Oversteegen.

\begin{defn} \cite{FO2002}\label{conj-equiv}
Two group chains $\{G_i\}_{i \geq 0}$ and $\{H_i\}_{i \geq 0}$ in $\fG$ are \emph{conjugate equivalent} if and only if there exists a collection $(g_i) \in G$, such that the group chains $\{g_iG_ig_i^{-1}\}_{i \geq 0}$ and $\{H_i\}_{i \geq 0}$ are equivalent. Here $g_iG_i = g_jG_i$ for all $i\geq 0$ and all $j \geq i$.
\end{defn}

The dynamical meaning of the equivalences in Definitions \ref{defn-greq} and \ref{conj-equiv} is given by the following theorem, which follows from results  in \cite{FO2002}; see also \cite[Section 3]{DHL2016a}.

\begin{thm}\label{equiv-rel-11}
Let $(G_\infty, G)$ and $(H_\infty, G)$ be inverse limit dynamical systems for group chains $\{G_i\}_{i \geq 0}$ and $\{H_i\}_{i \geq 0}$. Then we have:
\begin{enumerate}
\item \label{er-item1} The group chains $\{G_i\}_{i \geq 0}$ and $\{H_i\}_{i \geq 0}$ are equivalent if and only if there exists a homeomorphism $\tau: G_\infty \to H_\infty$ equivariant with respect to the $G$-actions on $G_\infty$ and $H_\infty$, and such that $\phi(eG_i) = (eH_i)$.
\item The group chains $\{G_i\}_{i \geq 0}$ and $\{H_i\}_{i \geq 0}$ are conjugate equivalent if and only if there exists a homeomorphism $\tau: G_\infty \to H_\infty$ equivariant with respect to the $G$-actions on $G_\infty$ and $H_\infty$.
\end{enumerate}
\end{thm}

That is, an equivalence of two group chains corresponds to the existence of a \emph{basepoint-preserving} conjugacy between their inverse limit systems, while a conjugate equivalence of two group chains corresponds to the existence of a conjugacy between their inverse limit systems, which need not preserve the basepoint.

Let $\{G_i\}_{i \geq 0} \in \fG$ with associated Cantor minimal system $(G_{\infty}, G, \Phi)$. 
An \emph{automorphism} of $(G_{\infty}, G, \Phi)$ is a homeomorphism $h \colon G_{\infty} \to G_{\infty}$ which commutes with the   $G$-action on $G_{\infty}$.  That is, for every $(g_{\ell}   G_{\ell}) \in G_{\infty}$ and   $g \in G$,    $g \cdot h(g_{\ell}   G_{\ell}) = h(g g_{\ell}   G_{\ell})$. 
Denote by ${Aut}(G_{\infty}, G, \Phi)$ the group of automorphisms of the action $(G_{\infty}, G, \Phi)$.  Note that ${Aut}(G_{\infty}, G, \Phi)$ is a topological group for  the compact-open topology on maps, and  is a closed subgroup of ${Homeo}(G_{\infty})$. 
 
 We introduce the definitions:
\begin{defn}\label{def-regular}
The Cantor minimal system  $(G_{\infty}, G, \Phi)$   is:
\begin{enumerate}
\item \emph{regular}  if the action of ${Aut}(G_{\infty}, G, \Phi)$ on $G_{\infty}$ is transitive;
\item  \emph{weakly regular} if the action of ${Aut}(G_{\infty}, G, \Phi)$ decomposes $G_{\infty}$ into a finite collection of  orbits;
\item \emph{irregular}  if the action of ${Aut}(G_{\infty}, G, \Phi)$ decomposes $G_{\infty}$ into an infinite collection of  orbits.
\end{enumerate}
\end{defn}

Let  $(G_\infty, G, \Phi)$ denote the Cantor minimal system associated to a group chain $\{G_i\}_{i \geq 0} \in \fG$.  Let $\fG(\Phi) \subset \fG$ denote the collection of all group chains in $\fG$ which are conjugate equivalent to $\{G_i\}_{i \geq 0}$.  
Then we have the following basic result, which follows from Theorem~\ref{equiv-rel-11}   above.  

\begin{thm}\label{thm-rephr}
 Let $\{G_i\}_{i \geq 0} \in \fG$ with associated Cantor minimal system $(G_{\infty}, G, \Phi)$.  Then $(G_{\infty}, G, \Phi)$ is:
\begin{enumerate}
\item \emph{regular} if all group chains in $\fG(\Phi)$ are equivalent;
\item  \emph{weakly regular} if $\fG(\Phi)$ contains a finite number of   classes of equivalent group chains;
\item \emph{irregular} if $\fG(\Phi)$ contains an infinite number of   classes of equivalent group chains.
\end{enumerate}
\end{thm}

We now consider the application of the above notions for the study of the homogeneity of weak solenoids.   Let $\cS_{\cP}$ be a weak solenoid, defined by a {presentation}    $\cP = \{ p_{\ell+1} \colon M_{\ell+1} \to M_{\ell} \mid \ell \geq 0\}$. Let $x \in \cS_{\cP}$ be a basepoint, and let $\{H_{i}^x\}_{i \geq 0}$ be the group chain in $H_0 = \pi_1(M_0 , x_0)$ defined by \eqref{eq-descendingchain}. 

 The following notion was introduced by McCord   in \cite{McCord1965}. 
  \begin{defn}\label{def-normal}
 A presentation $\cP$  is said to be   \emph{normal}  if  for each $\ell \geq 1$  the image    $H_{\ell}^x \subset H_0$ is a normal subgroup, and thus    each quotient  $X_{\ell}^x = H_{0}/H_{\ell}^x$ is      finite group, 
 \end{defn} 
 For example, if $H_0$ is an abelian group, then every group chain $\{H_{\ell}\}$ in $H_0$ is normal. 

For a normal  chain, the inverse limit \eqref{eq-Galoisfiber} is   a profinite group, which acts transitively on itself on the right. The   right action of $H_{\infty}$  on $H_{\infty}$,  commutes with the left action of $H_0$ on $H_{\infty}$,  and, moreover,  $H_{\infty} = {Aut}(H_{\infty}, H_0, \Phi)$. Thus, the automorphism group acts transitively on  $H_{\infty}$.
 McCord used this observation  in \cite{McCord1965} to show that the group ${Homeo}(\cS_{\cP})$ then acts transitively on $\cS_{\cP}$.

Rogers and Tollefson in \cite{RT1971b} gave an example of a weak solenoid for which the presentation is defined by a group chain which is not  normal, yet the inverse limit was still a profinite group, and so  the weak solenoid is homogeneous.  This example was the motivation for the work of  
Fokkink and Oversteegen in \cite{FO2002}, where they    gave a necessary and sufficient condition on the chain $\{H_i\}_{i \geq 0}$ for the weak solenoid to be homogeneous.
In this work, they showed the following result.
Denote by $N_{H_0}(H_i)$ the normalizer of $H_i$ in $H_0$, that is, $N_{H_0}(H_i) = \{h \in H_0 \mid h \, H_i \, h^{-1} = H_i\}$. 
Then we have the following result of Fokkink and Oversteegen \cite{FO2002}.
\begin{thm} \cite{Dyer2015} \label{thm-criteria}
Given a group chain  $\{H_i\}_{i \geq 0}$ with associated    equicontinuous Cantor minimal system  $(H_{\infty}, H_0, \Phi_0)$, then:
\begin{enumerate}
\item   $(H_{\infty}, H_0, \Phi_0)$ is regular if and only if there exists $\{N_i\}_{i\geq 0} \in \fG(\Phi_0)$ such that $N_i$ is a normal subgroup of $H_0$ for each $i \geq 0$. 
\item   $(H_{\infty}, H_0, \Phi_0)$ is weakly regular if and only if there exists $\{H_i'\}_{i\geq 0} \in \fG(\Phi_0)$ and an $n>0$ such that $H_i' \subset H_n \subseteq N_{H_0}(H_i')$ for all $i \geq n$. 
\end{enumerate}
\end{thm}

It follows that the class of weakly regular equicontinuous systems in Definition \ref{def-regular} is precisely the class of dynamical systems on the fibres of weakly normal solenoids in \cite{FO2002}. Note that the definition of equivalence of group chains in \cite{FO2002} is equivalent to ours, even if it is formulated slightly differently.
 Theorem~25 of \cite{FO2002} gave a criterion for when a weak solenoid is homogeneous, which combined with Theorem~\ref{thm-criteria} then yields:  
 \begin{thm}\label{thm-fokking}
 Let $\cS_{\cP}$ be a weak solenoid, defined by a {presentation}    $\cP$ with associated group chain $\{H_i\}_{i\geq 0}$. Then  $\cS_{\cP}$ is homogeneous if and only if  $\{H_i\}_{i\geq 0}$ is weakly regular. 
 \end{thm}

 The thesis \cite{Dyer2015} and the paper \cite{DHL2016a} considered an algebraic invariant for group chains, called the \emph{discriminant}. We do not include the definition here, but note that the study of the relation between the discriminant invariant and the homogeneity of a weak solenoid, leads to the following concept, which is the group chain version of  Definition~\ref{defn-normalcover} of the introduction. 
\begin{defn}\label{virtually-regular}
 The minimal equicontinuous action    $(H_{\infty}, H_0, \Phi_0)$   is \emph{virtually regular} if there is a   subgroup $H_0' \subset H_0$ of finite index   such that the   chain 
$\{H_i \cap H_0'\}_{i\geq 0}$ is regular as a chain in $H_0'$.
 \end{defn}
The virtually regular condition  in Definition~\ref{virtually-regular} appears to be similar to the    criteria in case (2) of Theorem~\ref{thm-criteria}. However, the work \cite{DHL2016a} gives several examples which show that the virtually regular property is independent from the weakly regular property.

\section{Germinal holonomy}\label{sec-holonomy}

In this section, we   introduce the notion of the \emph{germinal holonomy groups}  for an action $(\fX_0 , G_0, \Phi)$. The   definition is  analogous to that of the   germinal holonomy groups  of a foliated space.
 Given $x \in \fX_0$, let $G_x \subset G_0$ denote the isotropy subgroup of $x$.

\begin{defn}\label{def-holonomy} 
  Given  $g_1 , g_2 \in G_x$, we say $g_1$ and $g_2$ have the same  \emph{germinal  $\Phi$-holonomy} at $x$ if there exists an open set $U_x \subset \fX_0$ with $x \in U_x$, such that the restrictions $\Phi(g_1)|U_x$ and $\Phi(g_2)|U_x$  agree  on $U_x$.  In particular, we say that $g \in G_x$  has \emph{trivial  germinal $\Phi$-holonomy} at $x$ if there exists an open set $U_x \subset \fX$ with $x \in U_x$, such that the restriction $\Phi(g)|U_x$ is the trivial map.
\end{defn}

Note that the notion ``germinal $\Phi$-holonomy at $x$'' defines  an equivalence relation on the isotropy subgroup $G_x$. The reflexive and symmetric conditions are immediate, so we need only check the transitive condition. Suppose that  $g_1 , g_2, g_3 \in  G_x$ are such that $\Phi_0(g_1)$ and $\Phi_0(g_2)$ agree on the open set $U_x$, and $\Phi_0(g_2)$ and $\Phi_0(g_3)$ agree on the open set $V_x$,  then all three maps agree on the open set $U_x \cap V_x$ so that $g_1$ and $g_3$  have the same  germinal $\Phi_0$-holonomy at $x$. Thus, for each $x \in \fX$,   the  germinal $\Phi_0$-holonomy group at $x$ is defined, and denoted by ${\rm Germ}(\Phi_0 , x)$,    
and  there is surjective quotient map $G_x \to {\rm Germ}(\Phi_0 , x)$. 
 
We next recall a basic result of Epstein, Millet and Tischler \cite{EMT1977}.

\begin{thm} \label{thm-emt}
Let $(\fX_0 , G_0, \Phi)$ be a given action, and suppose that $\fX_0$ is a Baire space. 
Then the  union of all $x \in \fX_0$ such that  ${\rm Germ}(\Phi , x)$ is the trivial group forms a  $G_{\delta}$ subset of $\fX_0$. In particular, there exists at least one  $x \in \fX_0$ such that ${\rm Germ}(\Phi , x)$ is the trivial group.
\end{thm}
The main result in \cite{EMT1977} is stated in terms of the germinal holonomy groups of leaves of a foliation, but an inspection of the proof  shows that it  applies directly to an action $(\fX_0 , G_0, \Phi)$.

We   consider these notions for the special case of a Cantor minimal system $(G_{\infty}, G_0, \Phi)$ associated to a group chain $\cG = \{G_i\}_{i \geq 0} \in \fG$.
Let $K(\Phi) \subset G_0$ denote the kernel of $\Phi$, then $K(\Phi)$ is always a normal subgroup of $G_0$. The   action $\Phi$ is  \emph{effective} exactly when $K(\Phi)$ is the trivial subgroup.   
Consider the subgroup  $k(\cG) \subset G_0$ defined by 
\begin{equation}\label{eq-kernel}
k(\cG) = \bigcap_{i \geq 0} ~ G_i  ~ \subset G_0 ~ , 
\end{equation}
which is called the  \emph{kernel} of the chain $\cG$ in \cite{CortezPetite2008,Dyer2015,DHL2016a}.   
Note that  for any $g \in k(\cG)$ the   action $\Phi(g)$  fixes the     identity element $x_0 = (eG_i) \in G_{\infty}$, thus, $k(\cG)$ is the isotropy subgroup at $x_0$ for  $\Phi \colon G_0 \to {\rm Homeo}(G_{\infty})$. 
Note that if $\Phi(g)$ is the identity map for all $g \in k(\cG)$, that is, the action of every $g \in k(\cG)$ fixes every point in $G_{\infty}$, 
 then $k(\cG) = K(\Phi)$.

 Consider  $y \in G_{\infty}$, and let   $\cG_y = \{G_i^y\}_{i\geq 0}$ with $G_0^y = G$   be a group chain at $y$ associated to the action $(G_{\infty}, G_0, \Phi)$. 
Note that the kernel $k(\cG_y)$ of the group chain $\cG_y$ satisfies  $K(\Phi) \subseteq k(\cG_y)$.
 
 By Theorem~\ref{equiv-rel-11},  $\cG_y$ is    conjugate equivalent to $\cG$ as group chains in $G$.  If $y$ is in the orbit of the $G$-action  on $x_0$, then the kernels 
 $k(\cG)$ and $k(\cG_y)$ are  conjugate subgroups of $G$.   
 However, if $x_0$ and $y$ are contained in  distinct orbits of the $G$-action, then the groups $k(\cG)$ and $k(\cG_y)$ need not be conjugate.
There are examples of Cantor minimal systems     constructed   in \cite{DHL2016a,FO2002} from group chains,  such that $k(\cG)$ is not trivial, but there exists $y \in G_{\infty}$ such that $k(\cG_y)$ is   trivial. 

Note that if there exists $y, z \in G_{\infty}$ such that the kernel groups $k(\cG_y)$ and $k(\cG_z)$ are distinct, then the 
the action of ${Aut}(G_{\infty}, G_0, \Phi)$ on $G_{\infty}$ is not transitive, so the study of the family of subgroups 
$\{k(\cG_y) \mid y \in  G_{\infty}\}$ is a fundamental problem for the study of the dynamics of group chains.

We next explore a consequence of  
Theorem~\ref{thm-emt} for an action $(G_{\infty}, G_0, \Phi)$ associated to a group chain $\cG = \{G_i\}_{i \geq 0}$.  The theorem  implies that there exists $y \in G_{\infty}$ such that the germinal holonomy group ${\rm Germ}(\Phi , y)$ is trivial. Fix such a choice of $y$, and let  $\cG_y = \{G_i^y\}_{i\geq 0}$ with $G_0^y = G$   be a group chain at $y$ associated to the action $(G_{\infty}, G_0, \Phi)$, and let $k(\cG_y)$ be its kernel.  Then  $k(\cG_y)$ is the isotropy group at $y = (e G_i^y)$ of the action of $G$ on the inverse limit space $G_{\infty}^y$ associated with $\cG^y$.  
We   have  the following observation:
\begin{lemma}\label{lem-normalized}
For   $g \in k(\cG_y)$, there exists $i(g) \geq 0$ such that  $h^{-1} g h \in k(\cG_y)$ for all $h \in G^y_{i(g)}$.
\end{lemma}
\proof
  Let  $g \in k(\cG_y)$ with $\Phi(g) \ne Id$. 
Then $\Phi(g)$ is a homeomorphism of $G_{\infty}^y$ which fixes the basepoint $y = (eG^y_i)$.
As $y$ has trivial germinal holonomy, there exists an open neighborhood $U_y \subset G_{\infty}^y$ with $y \in U_i$ and $\Phi(g)$ is the identity map when restricted to $U_y$.

A point $z \in G_{\infty}^y$ is given by a sequence 
$$ z = (g_i G_i^y)  =  (G_0, g_1 G_1^y, g_2 G_2^y,\ldots) ~ {\rm where}~ g_j G_i^y = g_i G_i^y \textrm{ for all }j \geq i \ .$$
Define a descending chain $U_{\ell}$ of open neighborhoods of $(eG_i^y)$ as follows, for $\ell \geq 0$:
\begin{equation}\label{eq-opennbhds}
U_{\ell} = \{(G_0, g_1 G_1^y, g_2 G_2^y,\ldots) \mid g_j G_i^y = g_i G_i^y \textrm{ for all }j \geq i ~, ~ g_i = e ~{\rm for} ~ 1 \leq i \leq \ell  \} .
\end{equation}
Choose $i(g)$ sufficiently large so that $U_{i(g)} \subset U_y$.

By assumption, $k(\cG_y) \subset G_i^y$ for all $i \geq 0$. Thus, as the multiplication by $g$ is ``coordinate-wise''  on the sequences in  \eqref{eq-opennbhds}, for $g \in k(\cG_y)$  the action of $\Phi(g)$   satisfies $\Phi(g)(U_{\ell}) = U_{\ell}$ for all $\ell \geq 0$. 

The assumption that the action of $\Phi(g)$ is the identity on $U_y$ implies the same for $U_{i(g)}$, thus the restriction 
$\Phi(g) \colon U_{\ell} \to U_{\ell}$ is the identity map for   $\ell \geq i(g)$. That is, 
  for each $(h_i G_i^y) \in U_{i(g)}$ we have $g h_i G_i^y = h_i G_i^y$ for all $i\geq i(g)$.

Let $h \in G_{i(g)}^y$. Then define $(h_i G_i^y) \in U_{i(g)}$ by setting $h_i = e$ for $1 \leq i \leq  i(g)$, and $h_i = h$ for $i > i(g)$.  Then we have  $g h G_i^y = h G_i^y$ for all $i> i(g)$, so that 
  $h^{-1} g h G_i^y =  G_i^y$ for all $i> i(g)$.  Thus,  $h^{-1} g h \in G_i^y $ for all $i > i(g)$, hence  $h^{-1} g h \in k(\cG_y)$ as was to be shown.
\endproof

Now suppose that $k(\cG_y)$  is a finitely generated subgroup of $G$. Choose    generators $\{g_1, \ldots , g_d\}$ for  $k(\cG_y)$. Then by Lemma~\ref{lem-normalized}, for  each $1 \leq \ell \leq d$ there exists an index $i_{\ell}$ such that $h^{-1} g_{\ell} h \in k(\cG_y)$ for all $h \in G_{i_{\ell}}$. Set $i_0 = \max \{i_1 , \ldots , i_d\}$, then we have that $h^{-1} g_{\ell} h \in k(\cG_y)$ for all $h \in G_{i_0}$ and each $1 \leq \ell \leq d$. It follows that  $k(\cG_y)$ is a normal subgroup of $G_{i_0}^y$. We have shown:
\begin{lemma}\label{lem-normalized2}
Assume that   $k(\cG_y)$ is a finitely generated subgroup of $G_0$, then there exists $i_0 \geq 0$ such that  $k(\cG_y)$ is a normal subgroup of  $G_{i_0}^y$.
\end{lemma}

\section{Virtually nilpotent solenoids}\label{sec-nilpotent}

In this section, we give the definitions of nilpotent and virtually nilpotent weak solenoids, and show that these notions are homeomorphism invariants of weak solenoids.

Let $\cS_{\cP}$ be a weak solenoid with presentation $\cP$, and 
recall that  $\Pi_0  \colon \cS_{\cP} \to M_0$ is the projection onto the base of the solenoid.
Choose a basepoint $x_0 \in M_0$ then the fiber over $x_0$ is denoted by $\fX_0 = \Pi_0^{-1}(x_0)$, 
and the fundamental group of $M_0$ is denoted by  $H_0 = \pi_1(M_0, x_0)$.
The induced action on the fiber $\fX_0$ is given by the representation  $\Phi \colon H_0 \to {\rm Homeo}(\fX_0)$.

 \begin{defn}\label{def-nilpotent}
Let  $\cS_{\cP}$  be a weak solenoid with presentation $\cP$. We say that $\cS_{\cP}$   is \emph{nilpotent}, respectively \emph{abelian},  if the image 
$\Phi(H_0) \subset {\rm Homeo}(\fX_0)$ is   nilpotent, respectively is   abelian.
\end{defn}

We next introduce the virtual versions of these conditions, 
which are analogous to the notion of virtually regular in Definition~\ref{virtually-regular}.
 Let $H_0' \subset H_0$ be a subgroup of  finite index. Let    $p_0' \colon M_0' \to M_0$  be the finite  covering map associated to the subgroup $H_0' \subset H_0 = \pi_1(M_0 , x_0)$. Then for   $x_0' \in M_0'$ such that $p_0'(x_0') = x_0$,   the induced map on fundamental groups    $(p_0')_{\#} \colon G_0' = \pi_1(M_0' , x_0') \to H_0$ has image $H_0'$ by construction.
Let $\Pi_0' \colon \cS_{\cP}' \to M_0'$ be the weak solenoid which is the pull-back of $\cS_{\cP}$ via the map $p_0'$ of their base manifolds. 
Then we have an induced map $P_0' \colon \cS_{\cP}' \to \cS_{\cP}$ which satisfies the commutative diagram:

\begin{picture}(300,76)
\put(250,60){$\cS_{\cP}$}
\put(250,20){$M_0$}
\put(170,60){$\cS_{\cP}'$}
\put(170,20){$M_0'$}

\put(176,50){\vector(0,-1){17}}
\put(160,40){$\Pi_0'$}

\put(256,50){\vector(0,-1){17}}
\put(262,40){$\Pi_0$}

\put(210,67){$P_0'$}
\put(210,27){$p_0'$}
\put(200,22){$\vector(1,0){30}$}
\put(200,62){$\vector(1,0){30}$}
\end{picture}

The assumption that $p_0'$ is a finite covering map implies that the fibers of the induced map $P_0'$ have the same finite cardinality, hence $P_0' \colon \cS_{\cP}' \to \cS_{\cP}$ is a finite covering map of solenoids.   Let $\fX_0' = (\Pi_0')^{-1}(x_0')$ denote the fiber over $x_0'$, then the image   $P_0'(\fX_0') \subset \fX_0$ is a clopen subset. Let $(\fX_0', G_0', \Phi')$ be the induced action for $\cS_{\cP}'$.

\begin{defn}\label{def-virtuallynilpotent}
Let  $\cS_{\cP}$  be a weak solenoid. Then $\cS_{\cP}$   is \emph{virtually nilpotent}, respectively \emph{virtually abelian},  if there exists a subgroup $H_0' \subset H_0$   of  finite index such that the induced weak solenoid $\cS_{\cP}'$ is nilpotent, respectively is abelian. That is, the image $\Phi'(G_0') \subset {\rm Homeo}(\fX_0')$ is   nilpotent, respectively is   abelian. 
\end{defn}

 The following result shows the these virtual properties are homeomorphism invariants of weak solenoids. The proof omits some details that can be found in  the works \cite{ClarkHurder2013,CHL2013c}.
\begin{prop}\label{prop-virtuallynilpotent}
Suppose that $\cS_{\cP}$ and $\cS_{\cQ}$ are homeomorphic weak solenoids, then $\cS_{\cP}$ is virtually nilpotent, respectively virtually abelian,  if and only if $\cS_{\cQ}$ is as well. 
\end{prop}
\proof
Let $\cS_{\cP}$ and $\cS_{\cQ}$ be homeomorphic weak solenoids, defined by presentations $\cP$ and $\cQ$, respectively, with base manifolds $M_0$ and $N_0$ of the same dimension.  Choose basepoints $x_0 \in M_0$ and $y_0 \in N_0$, and let $\fX_{\cP} \subset \cS_{\cP}$ be the fiber over $x_0$, and $\fX_{\cQ} \subset \cS_{\cQ}$ the fiber over $y_0$. Let $G_0 = \pi_1(M_0 , x_0)$ and $H_0 = \pi_1(N_0, y_0)$. Let  $(\fX_{\cP} , G_0 , \Phi_{\cP})$ denote the monodromy action for $\cS_{\cP}$ and 
$(\fX_{\cQ} , H_0 , \Phi_{\cQ})$   the monodromy action for $\cS_{\cQ}$.

The actions  $(\fX_{\cP} , G_0 , \Phi_{\cP})$ and $(\fX_{\cQ} , H_0 , \Phi_{\cQ})$ are both minimal, so by \cite[Theorem~4.8]{CHL2013c} they are return equivalent. By definition, this means that there exists clopen subsets $U \subset \fX_{\cP}$ and $V \subset \fX_{\cQ}$ and a homeomorphism $\phi \colon U \to V$ which conjugates the induced pseudogroup action on $U$ to the induced pseudogroup action on $V$. Using the results of \cite[Sections~5,6]{CHL2013c}, the assumption that the action $(\fX_{\cP} , G_0 , \Phi_{\cP})$ is    equicontinuous   implies that there exists a clopen subset $U' \subset U$ so that for $V' = \phi(U') \subset V$, the sets 
\begin{eqnarray*}
G_{U'} & = & \{ g \in G_0 \mid \Phi_{\cP}(g)(U') = U'\} \subset G_0\\
H_{V'} & = & \{ h \in H_0 \mid \Phi_{\cQ}(h)(V') = V'\}  \subset H_0
\end{eqnarray*}
are subgroups of finite index. Then the image $\Phi_{\cP}(G_{U'}) \subset {\rm Homeo}(\fX_{\cP})$ maps the set $U'$ to itself. 
Similarly, the image $\Phi_{\cQ}(G_{V'}) \subset {\rm Homeo}(\fX_{\cQ})$ maps the set $V'$ to itself. 

Now, suppose that $\cS_{\cP}$ is virtually nilpotent, then there exists a subgroup of finite index $G_0' \subset G_0$ such that the image $\Phi_{\cP}(G_0') \subset {\rm Homeo}(\fX_{\cP})$ is a nilpotent group.  Then $G_0' \cap G_{U'}$ is also a subgroup of finite index in $G_0$, and its image $\G = \Phi_{\cP}(G_0' \cap G_{U'}) \subset {\rm Homeo}(\fX_{\cP})$   maps the set $U'$ to itself. 
As $\Phi_{\cP}(G_0')$ is a nilpotent group, every subgroup of it is also nilpotent, and so $\G$ is a nilpotent group.  

The conjugation map $\phi \colon U' \to V'$ then induces a subgroup $\phi^*(\G) \subset {\rm Homeo}(V')$ which is also nilpotent. As $G_0' \cap G_{U'}$ is a subgroup of finite index in $G_{U'}$, the image $\phi^*(\G)$ is a subgroup of finite index in $\Phi_{\cQ}(H_{V'}) \subset {\rm Homeo}(\fX_{\cQ})$. It follows that  $\cS_{\cQ}$ is virtually nilpotent. 

The reverse case, where we assume that $\cS_{\cQ}$ is virtually nilpotent, follows in the same way. The case of virtually abelian weak solenoids follows in exactly the same way.
\endproof

 \vfill
 \eject

\section{Proofs of main theorems}\label{sec-proofs}

We first give the proof of Theorem~\ref{thm-main0}. We assume that   $\fM$ is an equicontinuous matchbox manifold,  all of whose leaves have finite type and      polynomial growth.    Choose a homeomorphism $\fM \cong \cS_{\cP}$ with a weak solenoid, and choose a basepoint $x_0 \in M_0$ in the base manifold of the presentation $\cP$ and set $H_0 = \pi_1(M_0, x_0)$ and $\fX_0 = \Pi_0^{-1}(x_0)$.  Let $(\fX_0, H_0, \Phi)$ be the associated Cantor system.
 Then by Theorem~\ref{thm-emt}  there exists  $x \in \fX_0$  so   such that ${\rm Germ}(\Phi , x)$ is the trivial group. 
Let $L_x$ denote the leaf of the foliation $\F_{\cP}$ on $\cS_{\cP}$ containing $x$, then $L_x$ has no holonomy. 

By  Theorem~\ref{thm-coarseinvariance}, the assumption that all leaves of the foliation $\F$ on $\fM$ have polynomial growth, implies that the leaf $L_x$ also has polynomial growth type. That is, 
 there exists some integer $d_x \geq 1$ such that $L_x$ has polynomial growth type of degree at most $d_x$.

   Let  $\cH_x = \{H_i^x\}_{i \geq 0}$ be the group chain defined by \eqref{eq-descendingchain} for the basepoint $x$. Then by Theorem~\ref{thm-weaksuspensions}, there is a homeomorphism $L_x = \wtM_0/k(\cH_x)$. That is, $L_x$ is the covering space of $M_0$ associated to the subgroup $k(\cH_x) \subset H_0$. 
   Then by the Milnor-\v{S}varc Theorem~\ref{thm-Milnor}, the set $H_0/k(\cH_x)$ equipped with the word metric induced from $H_0$ also has polynomial growth type of degree at most $d_x$. 

The assumption that all leaves of $\fM$ have finite type implies that  $L_x$ has finite type, and thus the subgroup $k(\cH_x)$ is finitely generated.
 The group   ${\rm Germ}(\Phi , x)$ is trivial by the choice of $x$, so by Lemma~\ref{lem-normalized2}, there exists an index   $i_0 \geq 0$ such that $k(\cH_x)$ is a normal subgroup of $H_{i_0}$. 
 
 Choose a set of generators   $H_0^{(1)} = \{\gamma_1, \ldots, \gamma_{\nu}\} \subset H_0$ such that the restriction $H_0^{(1)} \cap H_{i_0}$ is also a generating set for $H_{i_0}$. 
The growth type of $H_0/k(\cH_x)$ is independent of the choice of generators for $H_0$, so the set  $H_0/k(\cH_x)$ again has polynomial growth type with respect to the generating set $H_0^{(1)}$. Then for the word metric on $H_0$ defined by the generating set $H_0^{(1)}$, which is used to define the metrics on quotient sets, the inclusion $H_{i_0}/k(\cH_x) \subset H_0/k(\cH_x)$ is an isometry. It then follows that $H_{i_0}/k(\cH_x)$ has polynomial growth type of degree at most $d_x$. 

The quotient set $H_{i_0}/k(\cH_x)$ has a   group structure, as $k(\cH_x)$ is a normal subgroup.  The generating set $H_0^{(1)} \cap H_{i_0}$ for $H_{i_0}$ descends to a generating set for the group $H_{i_0}/k(\cH_x)$, and thus this quotient group   has polynomial growth type of degree at most $d_x$. Thus, by  Theorem~\ref{thm-Gromov} of Gromov, there exists a nilpotent subgroup $\cN_x \subset H_{i_0}/k(\cH_x)$ of finite index.

Let $H_0' \subset H_{i_0}$ denote the preimage of $\cN_x$ under the quotient map $H_{i_0} \to H_{i_0}/k(\cH_x)$, so that 
$H_0'$ is a subgroup of finite index in $H_{i_0}$.
 As $H_{i_0} \subset H_0$ has finite index, we conclude that $H_0'$ also has finite index in $H_0$. 
 Let $p_0' \colon M_0' \to M_0$ be the finite   covering associated to   $H_0' \subset H_0$, and let  
$p_0' \colon \cS_{\cP}' \to \cS_{\cP}$ be the induced covering solenoid. 
 
 Choose a basepoint $x_0' \in M_0'$ such that $p_0'(x_0') = x_0$ and set $G_0 = \pi_1(M_0' , x_0')$.  Define the the group chain $\cG = \{G_i \}_{i \geq 0}$ where $G_i$ is the preimage in $G_0$ of the subgroup  $H_i \cap H_0'$.
Then    $\cS_{\cP}'$ is homeomorphic to the suspension space defined by the action of the fundamental group $G_0$ of $M_0'$    on the inverse limit space $G_{\infty}$ for the group chain $\{G_i \}_{i \geq 0}$. 

 Let  $\Phi' \colon H_0' \to {\rm Homeo}(G_{\infty})$ be the global holonomy map for the solenoid $\cS_{\cP}'$. 
 Let $K(\Phi') \subset G_0$ denote the kernel of the homomorphism $h_{\F}'$. 
Observe that by construction, the kernel  $k(\cG)$  is normal in $G_0$, and hence $k(\cG) = K(\Phi')$.  
 Thus, the image   $h_{\F}'(G_0)$ is isomorphic to   $H_0'/k(\cH_{x'}) \cong \cN_x$ which is a nilpotent group. 
This completes the proof of Theorem~\ref{thm-main0}.

We next give the proof of Theorem~\ref{thm-main1}, which uses the  the conclusions and notations of the proof of Theorem~\ref{thm-main0} above, and the algebraic results in Section~\ref{sec-growth}.

 Let $\fM$ be an equicontinuous matchbox manifold,   and suppose that all leaves of   $\F$    have   polynomial growth  at most $3$.  By Theorem~\ref{thm-main0}, we    obtain a nilpotent subgroup $\cN_x \subset H_{i_0}/k(\cH_x)$ of finite index, whose growth type in bounded above by the growth type of the leaf $L_x$. It is given that $L_x$ has polynomial  polynomial growth of degree most $3$, and thus $\cN_x$ is a nilpotent group of growth type at most $3$.  Thus, by   Proposition~\ref{prop-lowrank}, $\cN_x$ contains  a free abelian subgroup $\G_0' \subset \cN_x$ of finite index.
 
 If the subgroup $\G_0'$ is the trivial group, then $\cN_x$ is a finite group by  Lemma~\ref{lem-torsion}, which implies that $H_0$ is a finite group, which is impossible. Thus, we may assume that $\G_0'$ has rank at least $1$.
 
 Again using  the notation as in the proof of Theorem~\ref{thm-main0}, let $H_0' \subset H_{i_0}$ be the preimage of 
 $$\G_{0}' \subset \cN_x \subset H_{i_0}/k(\cH_x) \ .$$ 
 Then $H_0'$ has finite index in $H_{i_0}$ since $\G_0'$ has finite index in $H_{i_0}/k(\cH_x)$. As $H_{i_0}$ has finite index in $H_0$ we have that $H_0'$ has finite index in $H_0$.  
 Let $p_0' \colon M_0' \to M_0$ be the finite   covering associated to   $H_0' \subset H_0$, and let  
$p_0' \colon \cS_{\cP}' \to \cS_{\cP}$ be the induced covering solenoid. Choose a basepoint $x_0' \in M_0'$ such that $p_0'(x_0') = x_0$ and set $G_0 = \pi_1(M_0' , x_0')$.  Define the   group chain $\cG = \{G_i \}_{i \geq 0}$ where $G_i$ is the preimage in $G_0$ of the subgroup  $H_i \cap H_0'$.
Then    $\cS_{\cP}'$ is homeomorphic to the suspension space defined by the action of the fundamental group $G_0$ of $M_0'$    on the inverse limit space $G_{\infty}$ for the group chain $\{G_i \}_{i \geq 0}$. 
 The global holonomy of the weak solenoid $\cS_{\cP}'$ is conjugate to the image $h_{\F}'(H_0')$, which is   isomorphic  to $H_0'/k(\cH_{x'}) \cong \G_0'$ which is an abelian group. 

 By the proof of Theorem~\ref{thm-main0}, the group chain $\{G_i\}_{i \geq 0}$ associated to the weak solenoid $\cS_{\cP}'$ and given by $G_i = H_0' \cap H_i^x$ has kernel $k(\cH_{x'}) = K(\Phi')$, and the quotients $G_i/k(\cH_{x'})$ are abelian for $i \geq i_0$. Thus, the chain $\{G_i\}_{i \geq 0}$ is normal, so by Theorem~\ref{thm-criteria} and the result Theorem~\ref{thm-fokking} of Fokkink and Oversteegen \cite{FO2002}, we conclude that $\cS_{\cP}'$ is a homogeneous weak solenoid. 
 This completes the proof of Theorem~\ref{thm-main1}.

\begin{ex}{\rm 
We give an example of a solenoid which satisfies the hypotheses of Theorem~\ref{thm-main1}, and      is not homogeneous, but is virtually homogeneous.  The following construction can be viewed as an oriented version with trivial kernel of the example by   Rogers and Tollefson in \cite{RT1972}, as described   in \cite[page 3750]{FO2002}.

First, we construct the group   $H_0$.
Let $\mZ^2$ be the free abelian group of rank $2$, and we write $\mZ^2 = \{(a,b) \mid a,b \in \mZ\}$. Define an automorphism of $\mZ^2$ which acts on the generators by $A (a,b) = (-a, -b)$. That is, this automorphism is the rotation of the plane by the angle $\pi$, and so  $A$ has  order $2$. 

Define $H_0$ to be the split extension of $\mZ^2$ by $\mZ$, so we have a split exact sequence $\mZ^2 \to H_0 \to \mZ$, and so can write elements of $H_0$ as 
$H_0 = \{(a,b,c) \mid a,b,c \in \mZ\}$.  The quotient group $\mZ$ acts via conjugation on the normal subgroup $\mZ^2$ via the automorphism $A$. That is,
\begin{equation}\label{eq-mulitiplication}
(a,b,c) * (a',b',c') =  (a'',b'',c'') \quad {\rm where} ~ (a'', b'') = (a,b) + A^c(a',b') ~, ~ c'' = c + c' .
\end{equation}
The base manifold $M_0$ is a quotient of the product manifold $\mT^2 \times \mS^1$ by the action of an isometry $f$ of order $2$. Let $\phi_A \colon \mT^2 \to \mT^2$ be the action on $\mR^2/\mZ^2$ induced by the action $A (x,y) = (-x,-y)$ on $\mR^2$, which is orientation-preserving. Let $R_{\pi} \colon \mS^1 \to \mS^1$ be the rotation of the circle in the counterclockwise direction by the angle $\pi$, so $R_{\pi}$ also has order $2$. Then let $f = \phi_A \times R_{\pi}$ be the product action on $\mT^2 \times \mS^1$, and $M_0$ is the quotient manifold under this action. Note that $M_0$ is a $\mT^2$-bundle over $\mS^1$, where the monodromy action of the fundamental group of the base circle on the fundamental group of the fiber is  $A$.

Next,   define a group chain $\cH = \{H_n\}_{n \geq 0}$ in $\G$. Choose    primes $p, q > 3$ and set
\begin{equation}
H_n = \{(p^n a, q^n b , 3^n c) \mid a,b,c \in \mZ\} .
\end{equation}
The adjoint action $A$ preserves the subspace $\{(p^n a, q^n b , 0) \mid a,b \in \mZ\}$ so each $H_n$ is a subgroup.
Note that the kernel $k(\cH) = \{(0,0,0)\}$. 

Set  $H_0' = \{ (a,b,2c) \mid a,b,c \in \mZ \}$, which is a normal abelian subgroup of index $2$ in $H_0$. For $n\geq 0$, set: 
  $$G_n = H_n \cap H_0' = \{(p^n a, q^n b , 2 \cdot 3^n  c) \mid a,b,c \in \mZ\} .$$ 
The chain  $\{G_n\}_{n \geq 0}$ consists of free abelian groups, so is regular as a chain in $G_0 = H_0'$.  Hence, the   solenoid $\cS_{\cP}'$ formed from the group chain $\{G_n\}_{n \geq 0}$ is homogeneous, and thus the solenoid formed from the group chain $\{H_n\}_{n \geq 0}$ is virtually homogeneous. 
  On the other hand, we have:
\begin{lemma}
 The group chain $\cH = \{H_n\}_{n \geq 0}$   is irregular as a chain in $H_0$.
\end{lemma}
\proof
The  proof  follows the same ideas as used by   Fokkink and Oversteegen   in \cite[page 3750]{FO2002} to show that the Rogers and Tollefson example is irregular.

Fix $n \geq 1$ and let $\gamma = (0,0, 3^n) \in H_n$. Note that $\gamma^2 = (0,0,2 \cdot 3^n)$ is in the center of $\G$.

Let $e_1 = (1,0,0) \in H_0$ and $e_2 = (0,1,0) \in H_0$. Then $e_1^{-1} = (-1,0,0)$ and $e_2^{-1} = (0,-1,0)$. We calculate using \eqref{eq-mulitiplication} 
\begin{align*} e_1 \cdot (p^n a, q^n b , 3^n) \cdot e_1^{-1}   =    e_1 \cdot e_1 \cdot (p^n a, q^n b , 3^n)  =    (2 + p^n a, q^n b , 3^n)  ,\\
e_2 \cdot (p^n a, q^n b , 3^n) \cdot e_2^{-1}   =    e_2 \cdot e_2 \cdot (p^n a, q^n b , 3^n)  =    (p^n a, 2 + q^n b , 3^n)   . \end{align*}
Since $p, q$ are relatively prime to $2$,  the adjoint action $Ad(a e_1 + b e_2) \colon H_0 \to H_0$ does not fix the subgroup $H_n$ for  $0 < a < p^n$ and $0 < b < q^n$. 
However, the cosets $(a,b,0) \cdot H_n$ for $(a,b)$ of this form exhaust all left cosets of $H_n$ in $H_0$, so we have that 
 $N_{H_0}(H_n) = H_n$.  It then follows by the  criteria in part (2) of Theorem~\ref{thm-criteria} that the group chain $\cH$ cannot be weakly regular, and so the associated solenoid is not homogeneous
\endproof

 }
 \end{ex}

Finally, we give the proof of  Theorem~\ref{thm-main2}.   The conclusion in Theorem~\ref{thm-main1} that the solenoid $\cS_{\cP}$ is virtually homogeneous required that the induced group chain $\{G_i\}_{i \geq 0}$ is weakly regular, which follows trivially for the case of a group chain  in an abelian group $G_0$. However, Theorem~\ref{thm-main2} asserts that this conclusion fails without the abelian property on the chain, and so for the case where the growth type is polynomial of degree  greater than $3$.
  This follows as a consequence of a construction of an irregular group chain in a nilpotent group with polynomial growth rate 4. The construction of this example was   given    in the thesis of the first author \cite{Dyer2015}.

Let $\mH$ denote the continuous Heisenberg group, presented in the form $\mH = (\mR^3, *)$ with the group operation $*$ given by $(x,y,z)*(x',y',z')=(x+x',y+y',z+z'+xy')$. 
This operation is standard addition in the first two coordinates, with the added twist in the last coordinate.
The   subgroup   $\cH = (\mZ^3, *) \subset \mH$ is called the discrete Heisenberg group.   We think about $\cH$ as $\mZ^2 \times \mZ$, where  the  $\mZ^2 \times \{0\}$ factor is an abelian subgroup, and the $\{0\} \times \mZ$ factor has non-trivial commutator with the  first factor.  Then the coset space $M_0 = \mH/\cH$ is a compact $3$-manifold without boundary. Let $x_0 \in M_0$ be the coset of the identity $(0,0,0) \in \mH$, then $\cH$ is naturally identified with the fundamental group $\pi_1(M_0, x_0)$. 
We consider towers of coverings of $M_0$ defined by subgroup chains in $\cH$.

  The authors  Lightwood, {\c{S}}ahin and Ugarcovici considered in \cite{LSU2014} the   \emph{normal} subgroup chains in $\cH$, and gave a classification according to a scheme that was introduced in \cite{CortezPetite2008}. Our interest is in the chains in $\cH$ which are  irregular, and so not classified in  \cite{LSU2014}. The corresponding towers of coverings of $M_0$   then define weak solenoids, which need not be  homogeneous.

 We recall from \cite{LSU2014} the class of subgroups in $\cH$   which can be written in the form $\Gamma=\bA\mZ^2\times m\mZ$ where $\bA = \left( \begin{array}{cc} i & j \\ k & l\end{array}\right)$ is a 2-by-2 matrix with non-negative integer entries and $m>0$ is an integer. Then $\gamma \in \Gamma$ is of the form $\gamma = (ix+jy, kx+ly, mz)$ for some $x,y,z \in \mZ$. A straightforward computation gives the following lemma.

\begin{lemma}\label{lemma-Heissubgroup}
A set $\Gamma=\bA\mZ^2\times m\mZ$, where $\bA$ is a matrix with integer entries, and $m>0$ is an integer, is a subgroup if and only if $m$ divides both entries of one of the rows of $\bA$.
\end{lemma}

\begin{thm}\label{main-6}
Let $\cH$ be the Heisenberg group, and let $\bA_n=    \left( \begin{array}{cc} p^n & 0 \\ 0 & q^n\end{array}\right)$, where $p$ and $q$ are distinct primes. Then
the action represented by the group chain 
\begin{equation}\label{eq-defHn}
H_0 =\cH ~ , ~ \{H_n\}_{n \geq 1} = \{\bA_n \mZ^2 \times p^n \mZ\}_{n \geq 1}
\end{equation}
  is irregular, and has trivial kernel. 
\end{thm}

\proof
First consider an arbitrary subgroup $L = \bM \mZ^2 \times m\mZ$, where $m>0$ is an integer, and $\bM$ has integer entries. Let $h=(a,b,c) \in \cH$. Then for any $\gamma = (ix+jy, kx+ly, mz) \in L$, we compute $h*\gamma * h^{-1}$. First, we compute that $h^{-1}=(-a,-b,-c+ab)$. Then a straightforward computation shows that
\begin{equation}\label{eq-thirdentry}
  h*\gamma * h^{-1} = (ix+jy, kx+ly, mz+akx+aly-ixb-jyb) \ .
\end{equation}
Then $h*\gamma*h^{-1}$ is in $L$ only if $m$ divides $mz+akx+aly-ixb-jyb$. 
Then  set
$$ i = p^n,  j = 0,  k = 0,  l = q^n. $$
 Suppose that $\{H_n\}_{n \geq 0}$ is a weakly regular chain. Then  by Theorem \ref{thm-rephr}  there exists a subgroup of finite index $H \subset \cH$ and subscript $t \geq 0$, such that $H \supseteq H_t \supset H_{t+1} \supset \cdots$, and $\{H_n\}_{n \geq t}$ is equivalent to $\{g_nH_ng_n^{-1}\}_{n \geq t}$ for any choice of $(g_n) \subset N$ such that $g_sH_n = g_nH_n$ for any $s \geq n$. We can choose $H = H_s$, for any $s \geq t$.

We are now going to show that for every such $s$, there is an element $h \in H_s$ such that $\{H_n\}_{n\geq s}$ is not equivalent to  $\{h H_n h^{-1}\}_{n\geq s}$, and so $\{H_n\}_{n \geq 0}$ cannot represent a weakly regular action.

If $h = (p^sx',q^sy',p^s z') \in H_s$, and $\gamma = (p^nx,q^ny,p^nz) \in H_n$, then the third entry of $h * \gamma * h^{-1}$ in \eqref{eq-thirdentry} is equal to
    \begin{align}\label{eq-3entry} p^nz + p^sq^{n}x'y - q^sp^nxy' = p^s(p^{n-s} z + q^{n}x'y - q^s p^{n-s}xy'). \end{align}
 Let $x' = q$, so $h = (p^sq,q^sy',p^sz')$, and let $y = q$, so $\gamma = (p^nx, q^{n+1},p^nz)$. Then for any $n>s$ the element $h * \gamma * h^{-1}$ is never in $H_{s+1}$, since $p^{s+1}$ does not divide \eqref{eq-3entry}. So $\{H_n\}_{n \geq 0}$ is not weakly regular.
 
 It is clear from the definition \eqref{eq-defHn} that the kernel $\cap_{i\geq 1} \, H_n = \{0\}$.
\endproof

We now complete the proof of Theorem~\ref{thm-main2}. Associated to the chain   $\{H_n\}_{n \geq 0}$ is a sequence of compact quotient manifolds $M_n = \mH/H_n$ where $M_{n+1}$ is a finite covering of $M_n$ for $n \geq 0$. Let $\cP$ be the resulting presentation, and let  $\cS_{\cP}$ denote the weak solenoid defined by this inverse limit system. Then $\cS_{\cP}$ is not homogeneous by Theorem~\ref{thm-fokking}.

Note that each leaf of the foliation of $\cS_{\cP}$ is covered by the Heisenberg group  $\mH$,  and thus has growth type at most polynomial of degree $4$.   
To complete the proof  of  Theorem~\ref{thm-main2}, we show:

\begin{lemma}\label{lemma-notvh}
The solenoid $\cS_{\cP}$ is not virtually homogeneous.
\end{lemma}
\proof

Let $\cH' \subset \cH$ be a subgroup of finite index. We can assume that $\cH'$ is a normal subgroup by passing to its normal core.  Then consider the group chain $\{H_n'\}_{n \geq 1}$ in $\cH'$, where $H_n' = H_n \cap \cH'$. 
We show that  $\{H_n'\}_{n \geq 1}$ is not regular in $\cH'$.

 Consider the free abelian subgroup $\mZ^2 \subset \cH$ generated by the elements   $\{(1,0,0),(0,0,1)\}$. The intersection $\cH' \cap \mZ^2$ is then also free abelian, and so is generated by elements $\{(a,0,0),(0,0,0),(0,0,c)\}$ where $a,c$ are positive integers.
 We also have the free subgroup $\mZ \subset \cH$ generated by    the element    $\{(0,1,0)\}$, so its intersection with $\cH'$ is generated by $\{(0,b,0)\}$ for some positive integer $b$. Thus, every element of $\cH'$ has the form $(ax, by, cz)$ for $(x,y,z) \in \mZ^3$.

 \bigskip

Let $s >0$ be such that $\max\{ a,b,c\} < \min \{p^s,q^s\}$. Let 
$$A = lcm\{a,p^s\} ~ , ~ B  = lcm\{b,q^s\} ~ , ~ C = lcm\{c,p^s\} , $$ where $lcm$ denotes the least common multiple.   Then 
$$H_s' = \cH' \cap H_s = \{(Ax, By, Cz) \mid (x,y,z) \in \mZ^3\}  $$
and so for  $n \geq s$, we have 
$$H_n' = \cH' \cap H_n =  \{(A p^{n-s} x, B q^{n-s} y, C p^{n-s} z) \mid (x,y,z) \in \mZ^3\}  \  .$$

Note that if $\{H_n'\}_{n \geq s}$ is not regular in $H_s'$, then it is not regular in $\cH'$. We show that   $\{H_n'\}_{n \geq s}$ is not regular in $H_s'$, using an argument similar   to the one used to prove that $\{H_n\}_{n \geq 0}$ is not weakly regular in $\cH$.  It will follow  that    $\{H_n\}_{n \geq 0}$ is not virtually regular in $\cH$.
As before, see also Theorem \ref{thm-rephr}, we note that  $\{H_n'\}_{n \geq s}$ is regular in $H_s'$ if and only if for any sequence $(g_n) \in H_s'$ such that $g_n H_n' = g_m H_n'$ for all $m \geq n$ the chain $\{g_nH_n'g_n^{-1}\}_{n \geq s}$ is equivalent to $\{H_n'\}_{n \geq s}$.

If $h = (Aq,By',C z') \in H_s$, and $\gamma = (Ap^{n-s}x,Bq^{n-s}y,Cp^{n-s}z) \in H_n'$, then for the third entry of the expression \eqref{eq-thirdentry} we obtain 
$$ Cp^{n-s}z+ABq^{n-s}qy-ABp^{n-s}xy' \ . $$
Let $y = q$, and let $t \geq s$ be the power of $p$ in the prime decomposition of $AB$. Then $h (x,q,z)h^{-1}$ is not in $H_n'$ for $n \geq t$, hence for $n > t$ the chain $\{hH_n'h^{-1}\}_{n \geq s}$ is not equivalent to $\{H_n'\}_{n \geq s}$.
\endproof

 \section{Future research}\label{sec-future}

 There are a variety of questions about the geometry and classification of weak solenoids suggested by    Theorems~\ref{thm-main0}, \ref{thm-main1} and \ref{thm-main2} and their   proofs. In this section, we discuss three such areas of research: the role of the geometry of the base manifold $M_0$;  the algebraic classification of subgroup chains in the fundamental group $H_0 = \pi_1(M_0, x_0)$; and the relation between the discriminant invariant and   foliation theory for equicontinuous foliated spaces.

\subsection{$3$-manifolds} When the compact base manifold $M_0$ for a solenoid has dimension $2$, the only possible growth types are zero, linear, quadratic, or the fundamental group contains a free subgroup on two generators, and hence has   exponential growth. In the case where $M_0$ has dimension $3$, there are two more possibilities for the growth type of its fundamental group $H_0$, as described in Theorem~12.17  of  \cite{DK2013}:
 
 \begin{enumerate}
\item $H_0$ has growth type   polynomial of order at most $3$, and so Theorem~\ref{thm-main1} applies;
\item $H_0$ is a   nilpotent group with growth type polynomial of order $4$ and so  Theorem~\ref{thm-main1} applies;
\item $H_0$ is solvable with exponential growth type;
\item  $H_0$ contains a free subgroup on two generators, and hence has   exponential growth. 
\end{enumerate}
 An example of the groups of the third type are the Baumslag-Solitar groups
 $$H_0 = BS(1,p) = \langle a,b \mid aba^{-1} = b^p\rangle $$ 
which can be realized as the fundamental groups of compact $3$-manifolds via a geometric surgery construction. 

\begin{prob}\label{prob1}
Describe  the group chains for Baumslag-Solitar groups, or more generally, for solvable groups with exponential growth type.  Show that the discriminant invariant for all such group chains must be non-trivial.
\end{prob}

The fundamental group of a hyperbolic $3$-manifold  provides an  example of groups of the fourth type. Moreover these groups are always residually finite, so admit many group chains. That is, a hyperbolic $3$-manifold admits towers of finite coverings, not necessarily normal coverings, such that the inverse limit solenoid has a leaf with trivial fundamental group. 
Hyperbolic $3$-manifolds and their fundamental groups have been extensively studied, and their invariants have a central role in recent research in geometry. 
\begin{prob}\label{prob2}
Relate the discriminant invariant for a tower of finite coverings of a hyperbolic $3$-manifold with its traditional geometric invariants. 
\end{prob}

\subsection{Heisenberg manifolds} The examples constructed in the proof of  Theorem~\ref{thm-main2} are   Heisenberg manifolds of dimension $3$. For higher dimensions, it is also possible to construct analogous examples of Heisenberg manifolds which have fundamental groups with polynomial growth types, but there is no systematic understanding of which group chains in their fundamental groups have trivial or finite discriminant groups.   The following problem      lends itself nicely to an inductive approach by the growth rates of the groups.

\begin{prob}\label{prob3}
Characterize the group chains in general Heisenberg groups, such as was done for normal group chains in \cite{LSU2014}. Characterize the group chains with finite discriminant groups, and those for which the discriminant   is a Cantor group.
\end{prob}

\subsection{Foliated spaces}
 A weak solenoid is an equicontinuous matchbox manifold, and by Theorem~\ref{thm-equicontinuous} all such connected spaces are homeomorphic to a weak solenoid.  A matchbox manifold is a special case of a foliated space, as introduced by Moore and Schochet \cite{MS2006}, and the matchbox manifolds are characterized by the property that they are transversally totally disconnected. 
 The work of  {\'A}lvarez L{\'o}pez and  Moreira Galicia \cite{ALM2016} studied the properties of equicontinuous foliated spaces in general, from the point of view of their being the topological analog   of Riemannian foliations of smooth manifolds. In particular, their work developed a version of the Molino theory \cite{Molino1988} for the topological setting of foliated spaces. For the   case of matchbox manifolds, the work of the authors in \cite{DHL2016a} gives an alternate approach to constructing a Molino theory for these spaces, and in fact provides even stronger results for this special case as explained in \cite{DHL2016c}.

 The proof of Theorem~\ref{thm-main0} in Section~\ref{sec-nilpotent} implicitly uses  ideas of foliation theory, via Theorem~\ref{thm-emt}, and the study of the ends of leaves in weak solenoids in \cite{CFL2014} was also inspired by results from the theory of foliations of smooth manifolds. We formulate two questions, the first one, in a less general form, was asked  by Rogers and Tollefson \cite[Problem 2]{RT1971b}. 
 \begin{prob}\label{prob4}
Let $\fM$ be an equicontinuous matchbox manifold, hence a weak solenoid. Show that if $\fM$ is not a homogeneous space, then some leaf of its foliation $\F$ has non-trivial holonomy.
\end{prob}

Another related problem is based on the observation, that if a solenoid is homogeneous, then all leaves in this solenoid have either 1, or 2, or a Cantor set of ends. 
\begin{prob}\label{prob4}
Let $\fM$ be an equicontinuous matchbox manifold. Suppose that $\fM$ is virtually homogeneous, are their any restrictions imposed  on the end structures of its leaves? That is, is it possible to have a virtually homogeneous solenoid with a leaf for which the number of its ends is other than 1, 2 or a Cantor set?
\end{prob}

Note that Rogers and Tollefson example   is virtually homogeneous, and all of its leaves  have either 1 or 2  ends.   The Schori example has a leaf with 4 ends, but we do not know if it is virtually homogeneous.


\end{document}